\documentclass[leqno,12pt]{article}

\usepackage{amssymb}
\usepackage{euscript}
\usepackage[dvips]{graphicx}
\usepackage{flafter}
\usepackage{pstricks}
\usepackage[latin1]{inputenc}
\usepackage{amsmath}
\usepackage{amsfonts}
\usepackage{pstricks-add}

\usepackage{mathrsfs} 

\evensidemargin\oddsidemargin

\begin{document}

\baselineskip=18pt \setcounter{page}{1}

\renewcommand{\theequation}{\thesection.\arabic{equation}}
\newtheorem{theorem}{Theorem}[section]
\newtheorem{lemma}[theorem]{Lemma}
\newtheorem{proposition}[theorem]{Proposition}
\newtheorem{corollary}[theorem]{Corollary}
\newtheorem{remark}[theorem]{Remark}
\newtheorem{fact}[theorem]{Fact}
\newtheorem{problem}[theorem]{Problem}
\newtheorem{example}[theorem]{Example}

\newcommand{\eqnsection}{
\renewcommand{\theequation}{\thesection.\arabic{equation}}
    \makeatletter
    \csname  @addtoreset\endcsname{equation}{section}
    \makeatother}
\eqnsection

\def\r{{\mathbb R}}
\def\e{{\mathbb E}}
\def\p{{\mathbb P}}
\def\P{{\bf P}}
\def\E{{\bf E}}
\def\Q{{\bf Q}}
\def\z{{\mathbb Z}}
\def\N{{\mathbb N}}
\def\T{{\mathbb T}}
\def\G{{\mathbb G}}
\def\L{{\mathbb L}}

\def\deg{\chi}

\def\ee{\mathrm{e}}
\def\d{\, \mathrm{d}}
\def\S{\mathscr{S}}
\def\bs{{\tt bs}}



\vglue50pt

\centerline{\Large\bf Almost sure convergence for}

\bigskip

\centerline{\Large\bf stochastically biased random walks on trees}

\bigskip
\bigskip

\centerline{by}

\medskip

\centerline{Gabriel Faraud, $\;$Yueyun Hu $\;$and$\;$ Zhan Shi}

\medskip

\centerline{\it Universit\'e Paris XIII, Universit\'e Paris XIII \&
Universit\'e Paris VI}

\bigskip
\bigskip
\bigskip

{\leftskip=2truecm \rightskip=2truecm \baselineskip=15pt \small

\noindent{\slshape\bfseries Summary.} We are interested in the
biased random walk on a supercritical Galton--Watson tree in the
sense of Lyons~\cite{lyons90} and Lyons, Pemantle and Peres~\cite{lyons-pemantle-peres96},
and study a phenomenon of slow movement. In order to observe such a
slow movement, the bias needs to be random; the resulting random
walk is then a tree-valued random walk in random environment. We
investigate the recurrent case, and prove, under suitable general
integrability assumptions, that upon the system's non-extinction, the maximal displacement of the walk in the first $n$ steps, divided by $(\log n)^3$, converges almost surely to a known positive constant.

\bigskip


\bigskip


} 

\bigskip
\bigskip

\section{Introduction}
   \label{s:intro}

\subsection{Stochastically biased random walks on Galton--Watson trees}

$\phantom{aob}$Let $\T$ be a supercritical Galton--Watson tree
rooted at $\varnothing$. Two vertices $x$ and $y$ are said to be
connected, and denoted by $x\sim y$, if $x$ is either the parent or
a child of $y$. For a vertex $x\in \T$, we denote by $|x|$ the distance between $x$ and the root $\varnothing$, and $\varnothing=x_{0},x_{1},\dots,x_{|x|}$ the shortest path between the root and $x$. Let $\omega=(\omega(x), \, x\in \T \backslash \{
\varnothing\})$ be a sequence of vectors defined by $\omega(x) =
(\omega(x, \, y), \, y\sim x)$ such that $\omega(x, \, y)>0$,
$\forall y\sim x$ and that $\sum_{y\sim x} \omega(x, \, y)=1$.

Given the sequence $\omega$, we define a random walk $(X_n, \, n\ge
0)$ on $\T$ whose transition probabilities are
$$
P_\omega(X_{n+1} = y \, | \, X_n =x) = \omega(x, \, y).
$$

For each vertex $x\in \T \backslash \{ \varnothing\}$, we denote its
parent by ${\buildrel \leftarrow \over x}$, and its children by
$(x^{(1)}, \cdots, x^{(N(x))})$, where $N(x)$ stands for the number of children of $x$. Instead of looking at $\omega(x, \,
y)$ (for $y\sim x$ and $x\in \T$), it is often more convenient to
study $A(x) := (A_i(x), \, 1\le i\le N(x))$ defined by
\begin{equation}
    A_i(x) := {\omega(x, \, x^{(i)})\over \omega(x, \, {\buildrel
\leftarrow \over x})}, \qquad 1\le i\le N(x).
    \label{Ai}
\end{equation}

\medskip

{\noindent \bf Example } 
 {\bf (Biased random walk on a Galton--Watson tree).}
{\rm
 When $A_i(x)\equiv \lambda$, $\forall x$,
 $\forall i$ (where $0<\lambda<\infty$ is a
 constant), the random walk $(X_n)$ is the
 $\lambda$-biased random walk on $\T$ studied
 by Lyons~\cite{lyons90}, \cite{lyons92}, 
 Lyons, Pemantle and
 Peres~\cite{lyons-pemantle-peres-ergodic},
 \cite{lyons-pemantle-peres96}, Peres and
 Zeitouni~\cite{peres-zeitouni}, and Ben Arous et
 al.~\cite{benarous-fribergh-gantert-hammond}. More
 particularly, if
 $A_i(x)\equiv 1$, $\forall x$, $\forall i$, we get the simple
 random walk on $\T$.

 Ben Arous and Hammond~\cite{benarous-hammond}
 considered the case that $A_i(x)$ does not depend on
 $x$ nor on $i$, but can be random. They called the
 resulting walk $(X_n)$ randomly biased walk on $\T$,
 and proved that the walk is more regular in some sense
 than the biased random walk.\hfill$\Box$

} 


\medskip

We focus, in this paper, on a phenomenon of {\it slow movement} of
the walk in the recurrent case. In order to exhibit the slow
movement, the transition probabilities need to be random (which was
already the case in the aforementioned work of Ben Arous and
Hammond~\cite{benarous-hammond}): the resulting random walk $(X_n)$
is a so-called random walk in random environment. In dimension 1
(which, informally, corresponds to the case $N(x)=1$ for all $x$), a
celebrated theorem of Sinai~\cite{sinai} says that ${X_n\over (\log n)^2}$ converges in distribution to a non-degenerate law.

{F}rom now on, we assume $(A(x), \, x\in \T \backslash \{
\varnothing\})$ to be i.i.d., and write $A=(A_1, \cdots, A_N)$ for a
random vector having the distribution of any of $A(x)$.
Note that here $N$ itself is random and follows the law of reproduction of $\T$. We always use $\P$ to denote the probability with respect to the environment, and $\p := \P \otimes P_\omega$ the annealed measure. It is convenient to consider $(\omega, \T)$ as a marked tree, see (\ref{V}) below in terms of a branching random walk. As such, when we say ``for almost all environment $\omega$", we mean, in fact, for almost all $(\omega, \, \T)$.

Let us introduce the logarithmic moment-generating function of $A$:
$$
\psi(t) := \log \E \Big\{ \sum_{i=1}^N A_i^t\Big\} \in (-\infty, \,
\infty], \qquad t\ge 0.
$$

\noindent In particular, $\psi(0) = \log \E(N)$. We always assume
$\psi(0)>0$, so that the Galton--Watson tree is supercritical.
Furthermore, we assume that there exists some $\delta>0$ such that
\begin{equation}
    \psi(t) <\infty ,
    \; \; \forall \, t\in (-\delta, 1+\delta),
    \hbox{ \rm and that }
    \E(N^{1+\delta})<\infty.
    \label{H1}
\end{equation}

We first recall the following recurrence/transience
criterion:

\bigskip

\noindent {\bf Theorem A (Lyons and
Pemantle~\cite{lyons-pemantle}).} (i) If $\inf_{t\in [0, \, 1]}
\psi(t) \le 0$, then for almost all $\omega$, $(X_n)$ is recurrent.

(ii) If $\inf_{t\in [0, \, 1]} \psi(t) > 0$, then for almost all
$\omega$, $(X_n)$ is transient on the set of
non-extinction.

\bigskip

Theorem A was proved in \cite{lyons-pemantle} under the additional
condition that the distribution of $A_i$ does not depend on $i$;
this condition was removed in Faraud~\cite{faraud}. See also
Menshikov and Petritis~\cite{menshikov-petritis} for a proof of this
criterion (under the additional assumptions that $N>1$ is
deterministic and that the law of $A_i$ does not depend on $i$) via
Mandelbrot's multiplicative cascades.

The transient case (i.e., if $\inf_{t\in [0, \, 1]} \psi(t) > 0$)
has received much research attention recently (\cite{elie},
\cite{elie2}, \cite{benarous-fribergh-gantert-hammond},
\cite{benarous-hammond}).

If $\inf_{t\in [0, \, 1]} \psi(t) < 0$, the walk $(X_n)$ is positive
recurrent for almost all $\omega$; in this case, it is not hard (see
\cite{yztree}, under the additional assumptions that $N$ is
deterministic and that the law of $A_i$ does not depend on $i$) to
prove that ${1\over \log n}\max_{0\le k\le n}|X_k|$ converges almost
surely to a positive constant.

We assume $\inf_{t\in [0, \, 1]} \psi(t) = 0$ from now on. There are two different regimes in this case, depending on the sign of $\psi'(1)
= \ee^{-\psi(1)} \E\{ \sum_{i=1}^N A_i \log A_i\}$. If $\psi'(1)<0$
(and $\inf_{t\in [0, \, 1]} \psi(t) = 0$), then by defining $\kappa
:= \inf\{ t>1: \, \psi(t)=0\} \in (1, \, \infty]$ (with $\inf
\varnothing := \infty$) and $\kappa_1 := 1- {1\over \min\{ \kappa, \, 2\} }
\in (0, \, {1\over 2}]$, the order of magnitude of $|X_n|$ is,
loosely speaking, $n^{\kappa_1}$. [That $\kappa>1$ is a consequence of the convexity of $\psi$.] More precisely, as far as strong
convergence is concerned, Hu and Shi~\cite{yztree} proved (assuming
$N$ is deterministic and that the law of $A_i$ does not depend on
$i$) that $\max_{0\le k\le n}|X_k|=n^{\kappa_1 + o(1)}$, $\p$-almost surely. For (functional) weak convergence, Peres and
Zeitouni~\cite{peres-zeitouni} established a quenched functional
central limit theorem for biased random walk on $\T$ (corresponding
to the case $A_i = {1\over \E (N)}$, $\forall i$, and thus
$\kappa=\infty$; assuming moreover $N>1$ a.s.). The latter was
extended by Faraud~\cite{faraud} for walks satisfying $\kappa \in
(8, \, \infty]$ for the quenched case and $\kappa\in (5,\infty]$ for the annealed case. The problem of whether ${|X_n| \over n^{\kappa_1}}$
converges weakly (in either quenched or annealed setting) for the
whole region $\kappa \in (2, \infty] $ remains open, to the best of
our knowledge.

\begin{figure}
\centering

\scalebox{0.8} 
{
\begin{pspicture}(0,-1.91)(18.865833,1.91)
\rput(11.005834,-0.71){\psaxes[linewidth=0.04,labels=none,ticks=none,ticksize=0.10583333cm]{->}(0,0)(0,-1)(3,2)}
\rput(6.005833,-0.71){\psaxes[linewidth=0.04,labels=none,ticks=none,ticksize=0.10583333cm]{->}(0,0)(0,-1)(3,2)}
\rput(15.725834,-0.71){\psaxes[linewidth=0.04,labels=none,ticks=none,ticksize=0.10583333cm]{->}(0,0)(0,-1)(3,2)}
\usefont{T1}{ptm}{m}{n}
\rput(10.525833,-0.69){\footnotesize $0$}
\usefont{T1}{ptm}{m}{n}
\rput(15.345833,-0.69){\footnotesize $0$}
\usefont{T1}{ptm}{m}{n}
\rput(5.605833,1.53){\footnotesize $\psi(t)$}
\usefont{T1}{ptm}{m}{n}
\rput(10.385834,1.53){\footnotesize $\psi(t)$}
\usefont{T1}{ptm}{m}{n}
\rput(15.285833,1.53){\footnotesize $\psi(t)$}
\usefont{T1}{ptm}{m}{n}
\rput(8.995833,-0.99){\footnotesize t}
\usefont{T1}{ptm}{m}{n}
\rput(13.955833,-0.99){\footnotesize t}
\usefont{T1}{ptm}{m}{n}
\rput(18.775833,-0.99){\footnotesize t}
\usefont{T1}{ptm}{m}{n}
\rput(7.505833,1.73){\footnotesize $\psi'(1)<0$}
\usefont{T1}{ptm}{m}{n}
\rput(12.505833,1.71){\footnotesize $\psi'(1)=0$}
\usefont{T1}{ptm}{m}{n}
\rput(17.265833,1.71){\footnotesize $\psi'(1)>0$}
\usefont{T1}{ptm}{m}{n}
\rput(6.8458333,-1.07){\footnotesize 1}
\usefont{T1}{ptm}{m}{n}
\rput(11.845833,-1.07){\footnotesize 1}
\usefont{T1}{ptm}{m}{n}
\rput(16.605833,-1.07){\footnotesize 1}
\psbezier[linewidth=0.03](5.997918,0.19)(6.1658335,-0.11)(6.3258333,-0.33)(6.8858333,-0.73)(7.445833,-1.13)(8.485833,-1.35)(9.025833,-0.19)
\psline[linewidth=0.02cm,linestyle=dashed,dash=0.16cm 0.16cm](8.645833,-0.73)(8.645833,-1.63)
\usefont{T1}{ptm}{m}{n}
\rput(8.665833,-1.69){\footnotesize $\kappa$}
\psbezier[linewidth=0.03](11.005834,0.19945519)(10.965834,-0.2546154)(11.559343,-0.7013694)(11.845833,-0.70815384)(12.132323,-0.7149383)(13.025834,-0.77)(13.965833,0.57)
\psbezier[linewidth=0.03](15.725833,-0.05)(15.745834,-0.15)(15.785834,-0.67)(16.165833,-0.69)(16.545832,-0.71)(17.225834,-0.45)(18.585833,0.65)
\psline[linewidth=0.04cm](6.8858333,-0.61)(6.8858333,-0.87)
\psline[linewidth=0.04cm](11.885833,-0.59)(11.885833,-0.83)
\psline[linewidth=0.04cm](16.665833,-0.61)(16.665833,-0.83)
\usefont{T1}{ptm}{m}{n}
\rput(5.645833,-0.69){\footnotesize $0$}
\rput(0.98583335,-0.71){\psaxes[linewidth=0.04,labels=none,ticks=none,ticksize=0.10583333cm]{->}(0,0)(0,-1)(3,2)}
\usefont{T1}{ptm}{m}{n}
\rput(0.58583325,1.53){\footnotesize $\psi(t)$}
\usefont{T1}{ptm}{m}{n}
\rput(3.9758334,-0.99){\footnotesize t}
\usefont{T1}{ptm}{m}{n}
\rput(2.4858332,1.73){\footnotesize $\psi'(1)<0$}
\usefont{T1}{ptm}{m}{n}
\rput(1.8258333,-1.07){\footnotesize 1}
\psbezier[linewidth=0.03](0.96,0.17)(1.1369433,-0.13)(1.2327076,-0.27)(1.8388238,-0.69)(2.44494,-1.11)(3.2350614,-1.15)(4.18,-1.13)
\psline[linewidth=0.04cm](1.8658333,-0.59)(1.8658333,-0.85)
\usefont{T1}{ptm}{m}{n}
\rput(0.625833,-0.69){\footnotesize $0$}
\usefont{T1}{ptm}{m}{n}
\rput(2.46,1.37){\footnotesize $\kappa = \infty$}
\end{pspicture} 
}
\caption{The different possible shapes for $\psi$}
\end{figure}
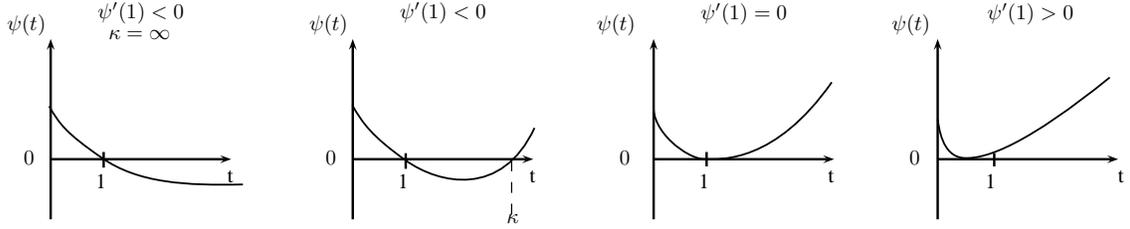

This paper is devoted to the study of the second situation:
$\psi'(1)\ge 0$ (and $\inf_{t\in [0, \, 1]} \psi(t) = 0$). It was
known that the walk is then extremely slow, at least under the
additional conditions that $N$ is deterministic and that the law of
$A_i$ does not depend on $i$: there exist constants $0<c_1\le
c_2<\infty$ such that (see \cite{yzslow})
\begin{equation}
    c_1 \le \liminf_{n\to \infty}\,
    {\max_{0\le k\le n} |X_k|\over (\log n)^3}
    \le \limsup_{n\to \infty}\,
    {\max_{0\le k\le n} |X_k|\over (\log n)^3}
    \le c_2, \qquad\hbox{\rm $\p$-a.s.}
    \label{yzslow}
\end{equation}

One of the main goals of this paper is to prove that almost sure
convergence holds in (\ref{yzslow}) if $\psi'(1)\ge 0$ (and
$\inf_{t\in [0, \, 1]} \psi(t) = 0$). The limiting constant in
(\ref{yzslow}), however, will have different natures depending on
whether $\psi'(1)=0$ or $\psi'(1)>0$. If $\psi'(1)\ge 0$ (and
$\inf_{t\in [0, \, 1]} \psi(t) = 0$), there exists $0<\theta\le 1$ such that
\begin{equation}
    \psi'(\theta)=0.
    \label{theta}
\end{equation}

\noindent The case $\theta=1$ corresponds to $\psi'(1)=0$; in this case, the condition $\inf_{t\in [0, \, 1]} \psi(t) = 0$ is equivalent to $\psi(1)=0$, and the walk is null recurrent. The case $\theta<1$, on the other hand, corresponds to $\psi'(1)>0$, and the walk is positive recurrent (see \cite{faraud}).

To give the limiting constant in (\ref{yzslow}), we define
\begin{equation}
    \alpha_\theta
    :=
    {3\pi^2\over 2} \, {1\over \theta}
    \, \E\Big( \sum_{i=1}^N A_i^\theta (\log A_i)^2 \Big) .
    \label{alpha}
\end{equation}

\noindent We write simply $\alpha$ for $\alpha_1$.

\begin{theorem}
\label{t:main}
 Assume $\psi(1) = \psi'(1)=0$. On the set of
 non-extinction,
 $$
 \lim_{n\to \infty}
 {\max_{0\le k\le n}|X_k| \over (\log n)^3}
 = {4\over \alpha} \, ,
 \qquad\hbox{\rm $\p$-a.s.,}
 $$
 where the constant $\alpha=\alpha_1$ is defined in $(\ref{alpha})$.

\end{theorem}

\begin{theorem}
\label{t:main2}
 Assume $\inf_{t\in [0, \, 1]} \psi(t) = 0$
 and $\psi'(1)>0$. On the set of
 non-extinction,
 $$
 \lim_{n\to \infty}
 {\max_{0\le k\le n}|X_k| \over (\log n)^3}
 = {1\over \alpha_\theta}   \, ,
 \qquad\hbox{\rm $\p$-a.s.,}
 $$
 where $\theta\in (0, \, 1)$ is as in
 $(\ref{theta})$ and $\alpha_\theta$ is defined in $(\ref{alpha})$.

\end{theorem}

\medskip

Theorem \ref{t:main2} is proved in Section \ref{s:proof-psi'(1)>0}, whereas the proofs of the upper and lower bounds in Theorem \ref{t:main} are in Sections \ref{s:ub} and \ref{s:lb}, respectively. In the next paragraph, we explain how it is that the $(\log
n)^3$ rate in Theorem \ref{t:main} arises, by means of an associated
branching random walk which plays the role of potential for our walk
$(X_n)$.

\subsection{Branching random walks and maxima along rays}

$\phantom{aob}$The influence of the random environment on the
behaviour of $(X_n)$ is best formulated in terms of an associated
potential process. To make the presentation easier, we artificially
add a special vertex, ${\buildrel \leftarrow \over \varnothing}$,
which is to be thought of as the parent of $\varnothing$. Since the
values of the transition probabilities at a finite number of vertices
have no influence on any of the results of the paper, we feel free
to modify the value of $\omega(\varnothing, \, \bullet)$, the
transition probability at $\varnothing$, in such a way that $(A_i(x),
\, 1\le i\le N(x))$, for $x\in \T$ (including $x= \varnothing$ now),
form an i.i.d.\ collection of random variables. Let
$\omega({\buildrel \leftarrow \over \varnothing}, \, \varnothing)=1$, and define consistently 
$$\omega(\varnothing,\, {\buildrel \leftarrow \over \varnothing}):=\frac{1}{1+\sum_{i=1}^{N(\varnothing)} A_{i}(\varnothing)}.$$

The potential process associated with the random environment is
defined by $V(\varnothing):=0$ and
\begin{equation}
    V(x) := -
    \sum_{y\in \, ]\!] \varnothing, \, x]\!]}
    \log \,
    {\omega({\buildrel \leftarrow \over y},
    \, y) \over
    \omega({\buildrel \leftarrow \over y}, \,
    {\buildrel \Leftarrow \over y})},
    \qquad x\in \T\backslash\{ \varnothing\},
    \label{V}
\end{equation}

\noindent where ${\buildrel \Leftarrow \over y}$ is the parent of
${\buildrel \leftarrow \over y}$, $[\![ \varnothing, \, x]\!]$
the set of vertices on the shortest path connecting $\varnothing$ to
$x$, and $\, ]\!] \varnothing, \, x]\!] := [\![ \varnothing, \, x]\!]
\backslash \{ \varnothing\}$.

Clearly, $(V(x), \, x\in \T)$ is a branching random walk, in the
usual sense of Biggins~\cite{biggins77}. It can be described as
follows: Initially, a single particle is located at $0$,
which is the ancestor of the system. At time 1, the ancestor dies,
giving birth to new particles who form the first
generation, and who are positioned according to the distribution of
$(-\log A_i(\varnothing), \, 1\le i\le N(\varnothing))$. At time 2, each
of the particles in the first generation dies, giving birth to new
particles that are positioned (with respect to their birth places)
according to the same distribution of $(-\log A_i(\varnothing), \,
1\le i\le N(\varnothing))$; these new particles form the second
generation. The system goes on according to the same mechanism. We
assume that for any $n$, each particle at generation $n$ produces
new particles independently of each other and of everything up to
the $n$-th generation. The positions of the particles in the $n$-th
generation are denoted by $(V(x), \, |x|=n)$.

The condition $\inf_{t\in [0, \, 1]} \psi(t)=0$ is equivalent to
$\inf_{t\in [0, \, 1]} \E(\sum_{|x|=1} \ee^{-t V(x)})=1$, whereas
$\psi'(1)\ge 0$ means $\E(\sum_{|x|=1} V(x) \ee^{- V(x)}) \le 0$.

In the recurrent case, there is a simple relationship between the
potential $(V(x), \, x\in \T)$ and the walk $(X_n)$. For any $k\ge
0$, let
$$
\tau_k := \inf\{ j\ge 1: \, |X_j| =k\}, \qquad \inf\varnothing :=
\infty.
$$

\noindent So $\tau_0$ is the first {\it return} time to the root if
the walk starts from $\varnothing$. It turns out that there exists
$0<c(\omega)<\infty$ possibly depending on the environment, such
that for any $n\ge 1$,
\begin{equation}
    \varrho_n:= P_\omega\{ \tau_n < \tau_0\}
    \ge {c(\omega)\over n} \exp\Big( - \min_{|x|=n}
    \overline{V}(x) \Big),
    \label{hitting-proba<=>potential}
\end{equation}

\noindent where, for any vertex $x$, we write
\begin{equation}
    \overline{V}(x) :=
    \max_{y\in \, ]\!] \varnothing, \, x]\!]} V(y).
    \label{Vbar}
\end{equation}

Inequality (\ref{hitting-proba<=>potential}) was proved in
\cite{yzslow} under the additional conditions that $N$ is
deterministic and that the law of $A_i$ does not depend on $i$.
Since the proof is simple, we reproduce it here: For any $x\in \T$,
let $T(x) := \inf  \{ i\ge 0: X_i =x \}$ be the first hitting time
of the walk at vertex $x$. By definition, for any $n\ge 1$, $\tau_n
= \min_{|x|=n} T(x)$, so that
\begin{equation}
    P_\omega\{ \tau_n < \tau_0\}
    \ge \max_{|x|=n} P_\omega \{ T(x) < \tau_0 \} .
    \label{beta>T}
\end{equation}

We fix a vertex $x$ with $|x|=n$. To compute $P_\omega \{ T(x) <
\tau_0 \}$, we define a random sequence $(\sigma_j)_{j\ge 0}$
(depending on $x$) by $\sigma_0:= 0$ and
$$
\sigma_j:= \inf \Big\{ k> \sigma_{j-1}: X_k \in [\! [\varnothing,\, x]
\! ] \backslash \{ X_{\sigma_{j-1}} \} \Big\}, \qquad j\ge 1.
$$

\noindent If the walk $(X_n)$ is recurrent, then $(\sigma_j)$ is
well-defined.

Let $Z_k := X_{\sigma_k}$, $k\ge 0$, which is the restriction of $(X_j)$ on the path $[\! [\varnothing,\, x] \! ]$. For $i\le n$, let $x_i$ be
the unique vertex in $[\![ \varnothing, \, x]\!]$ with $|x_i|=i$ (in
particular, $x_0=\varnothing$, $x_n=x$). Then for $1\le i<n$,
$$
P_\omega \Big\{ Z_{k+1}= x_{i+1} \, \Big| \, Z_k= x_i \Big\} =
{\omega(x_i, \, x_{i+1})\over \omega(x_i, \, x_{i+1}) + \omega(x_i,
\, x_{i-1})}
 = 1- P_\omega \Big\{ Z_{k+1}= x_{i-1} \,
\Big| \, Z_k= x_i \Big\} ,
$$

\noindent which yields
\begin{eqnarray}
    P_\omega \{ T(x) < \tau_0 \}
 &=&\omega(\varnothing, x_1) \, P_\omega \Big\{
    \hbox{$(Z_k)$ hits $x$ before
    hitting $\varnothing$} \, \Big| \, Z_0 = x_1 \Big\}
    \nonumber
    \\
 &=& {\omega(\varnothing, x_1) \, \ee^{V(x_1)} \over
    \sum_{z\in \, ]\! ] \varnothing,\, x] \! ]}
    \ee^{V(z)} } ,
    \label{zeitouni}
\end{eqnarray}

\noindent the second identity following from a general formula
(Zeitouni~\cite{zeitouni}, formula (2.1.4)) for the exit problem of
one-dimensional random walk in random environment. Since $\sum_{z\in \, ]\! ] \varnothing,\, x] \! ]}\ee^{V(z)} \ge \max_{z\in \, ]\! ] \varnothing,\, x] \! ]}\ee^{V(z)} = \ee^{\overline{V}(x)}$, this implies
\begin{equation}
    P_\omega \{ T(x) < \tau_0 \}
    \le
    \omega(\varnothing, x_1) \, \ee^{V(x_1)-\overline{V}(x)} .
    \label{zk-probatrans}
\end{equation}

\noindent Going back to (\ref{beta>T}), we immediately obtain
(\ref{hitting-proba<=>potential}) with $c(\omega) := \min_{|x|=1}
[\omega(\varnothing, x) \,\ee^{V(x)}] >0$.

The probability $\varrho_n$ is closely related to the maximal displacement of the branching random walk. The following simple observation was implicitly stated in \cite{yzslow} (pp.~1993--1996):

\begin{fact}
 \label{f:rhonxk}
 Assume $\inf_{t\in [0,1]} \psi(t)=0$ and $\psi'(1)\ge 0$. Let
 $0<c<\infty$ be a constant. Almost surely on the set of
 non-extinction,

 {\rm (i)} if $\varrho_n \ge \ee^{ - (c+o(1)) n^{1/3}}$ for all
 sufficiently large $n$, then
 $$
 \liminf_{n\to\infty} { 1\over (\log n)^3} \max_{0\le k\le n} |X_k|
 \ge {1\over c^3};
 $$

 {\rm (ii)} if $\varrho_n \le \ee^{ - (c+o(1)) n^{1/3}} $ for all
 sufficiently large $n$, then
 $$
 \limsup_{n\to\infty} { 1\over (\log n)^3} \max_{0\le k\le n} |X_k|
 \le {1\over c^3}.
 $$

\end{fact}

As such, an upper bound for $\min_{|x|=n}\overline{V}(x)$ yields, via inequality (\ref{hitting-proba<=>potential}), a lower bound for
$\varrho_n$, which, in turn, will lead to a lower bound for
the maximal displacement of the walk $(X_j)$.

\medskip

\begin{theorem}
\label{t:max-ray}
 Assume
 $\inf_{t\in [0, \, 1]}
 \E\{\sum_{|x|=1} \ee^{-t V(x)}\}=1$ and
 $\E\{\sum_{|x|=1} V(x) \ee^{- V(x)}\} \le 0$. Let
 $\theta\in (0, \, 1]$ be such that
 $\E\{\sum_{|x|=1} V(x) \ee^{- \theta V(x)}\} = 0$.
 We have, on the set of non-extinction,
 $$
 \lim_{n\to \infty} \,
 {1\over n^{1/3}} \min_{|x|=n} \overline{V}(x)
 = \Big( {3\pi^2 \sigma_\theta^2 \over 2}
 \Big)^{\! 1/3}
 \, , \qquad\hbox{\rm $\P$-a.s.,}
 $$
 where
 $$
 \sigma_\theta^2 :=  {1\over \theta}
 \E \Big\{ \sum_{|x|=1} V(x)^2
 \ee^{- \theta V(x)} \Big\} .
 $$

\end{theorem}

\medskip

We mention that Fang and Zeitouni \cite{FZ09} have independently obtained Theorem \ref{t:max-ray}, under the condition that $N$ is non-random and $A_i(\varnothing)$, for $1\le i\le N$, are i.i.d.

Comparing Theorem \ref{t:max-ray} with Theorems \ref{t:main} and \ref{t:main2}, we observe that (\ref{hitting-proba<=>potential}) is optimal in the case $\psi'(1)>0$ (or, equivalently, $\E\{ \sum_{|x|=1} V(x) \ee^{- V(x)}\} < 0$), but not in the case $\psi'(1)=0$ (or, equivalently, $\E\{ \sum_{|x|=1} V(x) \ee^{- V(x)} \} = 0$).

The proofs of the theorems are organized as follows.

$\bullet$ Section \ref{s:proof-max-ray}: Theorem \ref{t:max-ray}.

$\bullet$ Section \ref{s:proof-psi'(1)>0}: Theorem \ref{t:main2}.

$\bullet$ Section \ref{s:ub}: Theorem \ref{t:main}, upper bound.

$\bullet$ Section \ref{s:lb}: Theorem \ref{t:main}, lower bound. [This is the heart of the paper.]

\noindent Section \ref{s:spine} is devoted to a probability estimate for one-dimensional random walks, which will be exploited in the proofs of Theorems \ref{t:main} and \ref{t:main2} later on.

Throughout the paper, we use the convention $\sum_\varnothing := 0$, $\max_\varnothing := 0$ and $\min_\varnothing := \infty$. The letter $c$, with or without subscript, denotes a finite and positive constant, whose value may vary from line to line. Furthermore, $a_n \sim b_n$, $n\to \infty$, means $\lim_{n\to \infty} {a_n\over b_n} =1$.

\section{Proof of Theorem \ref{t:max-ray}}
\label{s:proof-max-ray}

$\phantom{aob}$Assume $\psi(1)=0$, i.e., $\E\{ \sum_{|x|=1} \ee^{- V(x)}\}=1$.

The condition $\E(N^{1+\delta})<\infty$ in (\ref{H1}) guarantees that $\P\{ N(x)<\infty, \; \forall x\} =1$ ($N(x)$ being the number of children of $x$). Recall that given a vertex $x\in \T$, $x_0:= \varnothing$, $x_1$, $\cdots$, $x_{|x|} := x$ are the vertices on $[\! [ \varnothing, \, x]\!]$ with $|x_i| =i$ for any $0\le i\le |x|$. The condition $\psi(1)=0$ yields that for any $n\ge 1$ and any measurable function $F: \, \r^n\times \r^n \to [0, \, \infty)$,
\begin{equation}
    \E\Big( \sum_{|x|=n} \ee^{-V(x)}
    F[V(x_i), \;
    N(x_{i-1}), \; 1\le i\le n] \Big)
    = \E \Big(
    F[S_i, \; \nu_{i-1}, \; 1\le i\le n]
    \Big),
    \label{change-proba2}
\end{equation}

\noindent where $(S_i-S_{i-1}, \; \nu_{i-1})$, for $i\ge 1$, are i.i.d.\ random vectors, whose common distribution is determined by
\begin{equation}
    \E[f(S_1, \; \nu_0)]
    =
    \E\Big( \sum_{|x|=1} \ee^{-V(x)}
    f ( V(x), \; N(\varnothing) ) \Big)
    =
    \E\Big( \sum_{i=1}^N A_i
    f ( -\log A_i, \; N ) \Big) ,
    \label{nu}
\end{equation}

\noindent for any measurable function $f: \, \r^2\to [0, \, \infty)$. Considering only the first argument, (\ref{change-proba2}) says that for any $n\ge 1$ and any measurable function $F: \, \r^n \to [0, \, \infty)$,
\begin{equation}
    \E\Big( \sum_{|x|=n} \ee^{-V(x)}
    F(V(x_i), \; 1\le i\le n) \Big)
    = \E [ F(S_i, \; 1\le i\le n) ],
    \label{change-proba}
\end{equation}

\noindent with the distribution of $S_1$ determined by
$$
\E(f(S_1) ) = \E\Big( \sum_{|x|=1} \ee^{-V(x)} f(V(x)) \Big) ,
$$

\noindent for any measurable function $f: \, \r\to [0, \, \infty)$. Formula (\ref{change-proba}) is well-known, and can be proved by means of a simple argument by induction in $n$. See, for example, Biggins and Kyprianou~\cite{biggins-kyprianou}. The proof of (\ref{change-proba2}) follows exactly from the same argument. In Section \ref{s:lb}, we will see an extension of (\ref{change-proba2}), which, in particular, gives a probabilistic interpretation of the new random walk $(S_i)$.

[The distribution of $S_1$ is well-defined upon the assumption $\psi(1)=0$. If furthermore $\psi'(1)=0$, then $\E(S_1)=0$; in words, $(S_n)$ is a mean-zero random walk, with $\sigma^2= \E(S_1^2)$.]

Formula (\ref{change-proba}) naturally leads to studying the
one-dimensional random walk $(S_n)$. However, we sometimes need to work in a slightly more general setting: For each $n\ge 1$, let $X_i^{(n)}$, $1\le i \le n$, be i.i.d.\ real-valued variables; define $S^{(n)}_0 :=0$ and $S^{(n)}_j := \sum_{i=1}^j X^{(n)}_i$ for $1\le j\le n$. Let $(a_n)$ be positive numbers such that $a_n \to\infty$ and ${a^2_n\over n} \to 0$, $n\to \infty$. Assume that there exists some $\eta>0$ and a  constant $\sigma^2>0$ such that, as $n \to\infty$,
\begin{equation}
    \E( X_1^{(n)})= o\Big( { a_n\over n} \Big), \qquad
    \sup_{n\ge 1}  \E ( |X_1^{(n)}|^{ 2+ \eta}) <\infty, \qquad
    \mbox{\rm Var}(X_1^{(n)}) \to \sigma^2.
    \label{hyp-array}
\end{equation}

\noindent The following estimate is essentially due to
Mogulskii~\cite{mogulskii}:

\begin{proposition}
 [A triangular version of Mogulskii~\cite{mogulskii}]
 \label{p:mogulskii-array}
 Assume (\ref{hyp-array}).
 Let $g_1< g_2$ be continuous functions on
 $[0, \, 1]$ with $g_1(0) < 0 < g_2(0)$.
 Consider the measurable event
 $$
 F_n:= \Big\{ g_1 ({i\over n}) \le
 {S^{(n)}_i\over a_n} \le g_2 ({i\over n}), \;
 \hbox{for}\; 1\le i\le n \Big\}.
 $$
 We have
 \begin{equation}
     \lim_{n\to \infty} {a_n^2\over n} \log
     \P(F_n) =
     - {\pi^2 \sigma^2\over 2}
     \int_0^1 {\d t\over [g_2(t)- g_1(t)]^2} .
     \label{mogulskii-array}
 \end{equation}
  Moreover, for any $b>0$,
 \begin{equation}
      \lim_{n\to \infty} {a_n^2\over n} \log
      \P\Big\{ F_n, \; {S^{(n)}_n\over a_n} \ge
     g_2(1)-b\Big\} =
      - {\pi^2 \sigma^2\over 2}
      \int_0^1 {\d t\over [g_2(t)- g_1(t)]^2} .
      \label{mogulskii-array2}
  \end{equation}
\end{proposition}

If the law of $X_1^{(n)}$ does not depend on $n$ (in which case we can even take $\eta=0$), Proposition \ref{p:mogulskii-array} is Mogulskii~\cite{mogulskii}'s theorem. For a detailed proof of Proposition \ref{p:mogulskii-array}, see \cite{GHS}.

A useful consequence of Proposition \ref{p:mogulskii-array} is as follows. Again, if the law of $X_1^{(n)}$ does not depend on $n$, we only need $X_1^{(n)}$ to have a finite second moment in order to have (\ref{hyp-array}).

\begin{corollary}
 \label{c:tech1}
 Assume that (\ref{hyp-array}) is satisfied with $a_n := n^{1/3}$.

 {\rm (i)} Let $f: [0, \, 1] \to (0, \infty)$ be a continuous function,
 and let $(f_n)$ be a sequence of continuous functions converging
 uniformly to $f$ on $[0, \, 1]$. Then for any $b >0$, when
 $n\to \infty$,
 $$
 \sup_{0\le u \le b \, n^{1/3}}
 \P\Big( u \ge S^{(n)}_i \ge  u - n^{1/3} f_n( {i\over n}), \;
 1\le i \le n\Big)
 =
 \ee ^{ -  {\pi^2  \sigma^2 \over 2}
 (1+o(1)) n^{1/3} \int_0^1 {\!\d t \over f^2(t)}}.
 $$

 {\rm (ii)} For any $b>a>0$, we have, as $n \to \infty$,
 $$
 \sum_{j=1}^n  \ee^{-b (n-j)^{1/3}}
 \P\Big( a n^{1/ 3} \ge S^{(n)}_i > a n^{1/ 3} - b (n-i)^{1/3}, \;
 \forall\, 1\le i \le j\Big)
 =
 \ee^{- \min\{ b, \, {3\pi^2 \sigma^2 \over 2 b^2}\} (1+o(1)) n^{1/3}}.
 $$

\end{corollary}

\noindent {\it Proof of Corollary \ref{c:tech1}.} We first prove (ii). Let $\varepsilon>0$. Define $k:=\lfloor {1\over \varepsilon} \rfloor$, $n_\ell:=\ell \lfloor \varepsilon n  \rfloor$ for $\ell =0$, $\cdots$, $k-1$ and $n_k:=n$. By (\ref{mogulskii-array}), the sum in (ii) is, for all large $n$ and some constant $c$,
\begin{eqnarray*}
 &\le&  \lfloor \varepsilon n  \rfloor \, \sum_{\ell=1}^k \ee^{- b( n- n_\ell)^{1/3}}\,
     \P\Big( a n^{1/ 3} \ge S^{(n)}_i > a n^{1/ 3} - b (n-i)^{1/3}, \;
     \forall\, 1\le i \le n_{\ell-1}\Big)
     \\
 &\le& \lfloor \varepsilon n  \rfloor \,  \sum_{\ell=1}^k \ee^{- b( n- n_\ell)^{1/3}}\,
    \ee^{ - ({3\pi^2\sigma^2 \over 2b^2} - \varepsilon)
    ( n^{1/3} - ( n- n_{\ell-1})^{1/3})}
    \\
 &\le&\ee^{- \min\{ b,\, {3\pi^2 \sigma^2 \over 2 b^2}\}
    (1-c\varepsilon) n^{1/3}} .
 \end{eqnarray*}

\noindent This proves the upper bound in (ii) as $\varepsilon$ can be arbitrarily small. The lower bound is easier: we only need to consider two terms: $j= \lfloor \varepsilon n  \rfloor$ and $j=n$, and apply again (\ref{mogulskii-array}).

The proof of (i) goes along similar lines by cutting the interval $\{0\le u \le b \, n^{1/3}\}$ into smaller intervals of length of order $\varepsilon n$ with small $\varepsilon>0$, using monotonicity and applying Proposition \ref{p:mogulskii-array}. The details are omitted. \hfill $\Box$

\bigskip

We now proceed to the proof of Theorem \ref{t:max-ray}: if $\inf_{t\in [0, \, 1]} \E\{ \sum_{|x|=1} \ee^{-t V(x)}\} =1$ and $\E\{ \sum_{|x|=1} V(x) \ee^{- V(x)} \} \le 0$, then on the set of non-extinction,
$$
\lim_{n\to \infty} \, {1\over n^{1/3}} \min_{|x|=n} \overline{V}(x)
=  \Big( {3\pi^2 \sigma_\theta^2 \over 2} \Big)^{\! 1/3} \, ,
\qquad\hbox{\rm $\P$-a.s.,}
$$

\noindent  where $\sigma_\theta^2 :=  {1\over \theta} \E \{ \sum_{|x|=1} V(x)^2 \ee^{- \theta V(x)} \}$ and $\theta\in (0, \, 1]$ is such that $\E \{\sum_{|x|=1} V(x) \ee^{- \theta V(x)}\} = 0$.

Without loss of generality, we can assume $\theta=1$. Indeed, if $0<\theta<1$, then by considering $\widetilde{V}(x) := \theta V(x)$, we have $\inf_{t\in [0, \, 1]} \E(\sum_{|x|=1} \ee^{-t \widetilde{V}(x)})=1$ and $\E(\sum_{|x|=1} \widetilde{V}(x) \ee^{- \widetilde{V}(x)}) =0$, so that by the case $\theta=1$, ${1\over n^{1/3}} \min_{|x|=n} \max_{y\in \, ]\!] \varnothing, \, x]\!]} \widetilde{V}(y) \to ( {3\pi^2 \widetilde{\sigma}^2 \over 2} )^{1/3}$ $\P$-almost surely on the set of non-extinction, where $\widetilde{\sigma}^2 := \E \{ \sum_{|x|=1} \widetilde{V}(x)^2 \ee^{-\widetilde{V}(x)} \}$.

So we only need to prove Theorem \ref{t:max-ray} in the case
$\theta=1$. In the rest of the section, we assume $\E(\sum_{|x|=1} \ee^{- V(x)})=1$ and $\E(\sum_{|x|=1} V(x) \ee^{- V(x)})=0$, and prove that, on the set of non-extinction,
\begin{equation}
    \lim_{n\to \infty} \,
    {1\over n^{1/3}} \min_{|x|=n} \overline{V}(x)
    = \Big( {3\pi^2 \sigma^2 \over 2} \Big)^{\! 1/3}
    \, , \qquad\hbox{\rm $\P$-a.s.},
    \label{max-ray}
\end{equation}

\noindent with $\sigma^2 := \sigma_1^2 = \E \{ \sum_{|x|=1} V(x)^2 \ee^{- V(x)} \}$. For the sake of clarity, we prove the upper and lower
bounds in distinct parts.

\bigskip

\noindent {\it Proof of (\ref{max-ray}): lower bound.} We assume $\E(\sum_{|x|=1} \ee^{- V(x)})=1$ and $\E(\sum_{|x|=1} V(x) \ee^{- V(x)})=0$.

Let $0<a < ({3 \pi^2 \sigma^2 \over 2})^{1/3}$ and $b := ({3 \pi^2 \sigma^2 \over 2})^{1/3}$. Let $n\ge 1$. For all $|x|=n$, let
$$
H_x:= \inf\{j \in [ 1, n] : V(x_j) \le a n^{1/3} - b (n-j)^{1/3} \}, \qquad \inf \varnothing :=\infty.
$$

\begin{figure}[h!]
\centering
\scalebox{1} 
{
\begin{pspicture}(0,-3.0725)(11.16,3.1125)
\rput(2.16,-0.0725){\psaxes[linewidth=0.04,labels=none,ticks=none,ticksize=0.10583333cm]{->}(0,0)(-1,-3)(9,3)}
\psline[linewidth=0.04](2.18,-0.0675)(2.64,-0.3925)(3.04,-0.5925)(3.44,-0.1925)(3.84,-0.3925)(4.24,-0.1925)(4.64,-0.7925)(5.04,-1.1925)(5.44,-0.7925)(5.84,0.0075)(6.44,-1.1925)(6.84,-0.7925)(7.04,0.2075)(7.36,0.7125)
\psline[linewidth=0.04cm,linestyle=dashed,dash=0.16cm 0.16cm](2.16,1.9275)(10.16,1.9275)
\psline[linewidth=0.04cm](10.16,0.1275)(10.16,-0.2725)
\usefont{T1}{ptm}{m}{n}
\rput(10.19,-0.4875){$n$}
\usefont{T1}{ptm}{m}{n}
\rput(1.62,2.9175){$V(x_j)$}
\usefont{T1}{ptm}{m}{n}
\rput(1.6,2.0175){$an^{1/3}$}
\psline[linewidth=0.04](7.36,0.7275)(7.76,-0.0725)(8.16,0.9275)(8.56,1.5275)(8.96,0.9275)(9.36,0.3275)
\usefont{T1}{ptm}{m}{n}
\rput(1.11,-1.0225){$(a-b)n^{1/3}$}
\psline[linewidth=0.04cm](1.96,-1.0725)(2.36,-1.0725)
\psline[linewidth=0.04cm](9.36,0.3275)(9.76,0.9275)
\psline[linewidth=0.04cm](9.76,0.9275)(10.16,1.3275)
\psline[linewidth=0.04cm,arrowsize=0.05291667cm 2.0,arrowlength=1.4,arrowinset=0.4]{<-}(4.66,-0.7425)(6.56,-2.0725)
\usefont{T1}{ptm}{m}{n}
\rput(6.45,-2.3675){$H_x$}
\psplot[linewidth=0.04cm,plotpoints=200]{2.1}{10.1}{2 8 x 2.1 sub sub 1 3 div exp 1.5 mul sub 0.09 sub}
\end{pspicture}
}

\caption{$H_{x}$}
\end{figure}

Assume there exists a vertex $x$ with $|x|=n$ such that $\overline{V}(x) \le a n^{1/3}$. Then $H_x \le n$; writing $j:= H_x$ and $y := x_j$, we have, for all $i<j$, $a n^{1/3} \ge V(y_i) > a n^{1/3} - b (n-i)^{1/3}$ and $V(y) \le a n^{1/3} - b (n-j)^{1/3}$. Therefore, by writing
$$
U_j:= \sum_{|y|=j} {\bf 1}_{
\{V(y) \le a n^{1/3} - b (n-j)^{1/3}, \, \,  a n^{1/3} \ge V(y_i)
> a n^{1/3} - b (n-i)^{1/3}, \, \forall i<j \}} ,
$$

\noindent we obtain:
$$
\P\Big( \min_{|x|=n} \overline{V}(x) \le a n^{1/3} \Big) \le \P \Big( \bigcup_{j=1}^n \{  U_j  \ge 1\}\Big) \le \sum_{j=1}^n \E ( U_j ).
$$

\noindent By (\ref{change-proba}), we have $\E(U_j) = \E[\ee^{S_j} {\bf 1}_{\{S_j \le a n^{1/3} - b (n-j)^{1/3}, \; a n^{1/3} \ge S_i > a n^{1/3} - b (n-i)^{1/3}, \; \forall i<j \}}]$. Hence
$$
\P\Big( \min_{|x|=n} \overline{V}(x) \le a n^{1/3}
\Big) \le \sum_{j=1}^n \ee^{a n^{1/3} - b (n-j)^{1/3}} \P\Big(a
n^{1/3} \ge S_i
> a n^{1/3} - b (n-i)^{1/3}, \, \forall i<j \Big).
$$

\noindent Applying Corollary \ref{c:tech1} (ii) and noting that $\min \{ b, \, {3\pi^2 \sigma^2 \over 2 b^2}\} = ({3\pi^2 \sigma^2 \over 2})^{1/3}$, we get that, for any $0< a < ({3 \pi^2 \sigma^2 \over 2})^{1/3}$,
\begin{equation}
    \limsup_{n \to\infty}
    {1\over n^{1/3}} \log
    \P\Big( \min_{|x|=n} \overline{V}(x) \le a n^{1/3} \Big)
    \le
    a - \Big( {3 \pi^2 \sigma^2 \over 2} \Big)^{1/3},
    \label{V<}
\end{equation}

\noindent which implies $\sum_n \P\{ \min_{|x|=n} \overline{V}(x) \le a n^{1/3} \}<\infty$. The lower bound in (\ref{max-ray}) follows from the Borel--Cantelli lemma, as $a$ can be as close to $({3 \pi^2 \sigma^2 \over 2})^{1/3}$ as possible.\hfill $\Box$

\bigskip

\noindent {\it Proof of (\ref{max-ray}): upper bound.} Assume
$\E\{ \sum_{|x|=1} \ee^{- V(x)}\}=1$ and $\E\{\sum_{|x|=1} V(x)
\ee^{- V(x)}\}=0$.

Let $n\ge 1$ and $b>a>\varepsilon>0$. The key step in the proof of the upper bound in (\ref{max-ray}) is the following estimate, which is a consequence of the Paley--Zygmund inequality (see \cite{GHS} for a proof): For any Borel sets $I_{i,n}\subset \r$, $1\le i\le n$, and any integer $r_n \ge 1$, we have
\begin{equation}
    \P\Big\{
    \exists |x|=n : \; V(x_i) \in I_{i,n} \, , \;
    \forall 1\le i\le n\Big \}
    \ge
    {\E [\ee^{S_n} \, {\bf 1}_{\{ S_i \in I_{i,n}\, , \;
    \nu_{i-1} \le r_n, \; \forall \, 1\le i\le n\} } ]
    \over 1+ (r_n-1)\sum_{j=1}^n h_{j,n}} ,
    \label{lemme-second-moment}
\end{equation}

\noindent where
$$
h_{j,n}:= \sup_{u\in I_{j,n}} \E \Big( \ee^{S_{n-j}} {\bf 1}_{ \{ S_i \in I_{i+j,n}-u, \; \forall 0\le i\le n-j\} } \Big) ,
$$

\noindent and $I_{i+j,n}-u := \{ v-u: \; v\in I_{i+j,n} \}$. [We recall that $(S_i-S_{i-1}, \, \nu_{i-1})$, $i\ge 1$, are i.i.d.\ random vectors (with $S_0:=0$) whose common distribution is given by $(\ref{nu})$.]

We choose $r_n := \lfloor \ee^{n^{1/4}}\rfloor$ and $I_{i, n}:=
[(a-\varepsilon) n^{1/3}-b(n-i)^{1/3}, \, a n^{1/3}]$. In particular, $\{ \exists |x|=n : \; V(x_i) \in I_{i,n} \, , \; \forall 1\le i\le n\} \subset \{ \min_{|x|=n} \overline{V} (x) \le an^{1/3}\}$. It follows from (\ref{lemme-second-moment}) that
$$
\P \Big( \min_{|x|=n} \overline{V}(x) \le a n^{1/3} \Big)
\ge
{\ee^{(a-\varepsilon) n^{1/3}} \, \P\{ S_i \in I_{i,n}\, , \;
\nu_{i-1} \le \ee^{n^{1/4}}, \; \forall \, 1\le i\le n\}
\over 1+ \ee^{n^{1/4}}\sum_{j=1}^n h_{j,n}} .
$$

Let $X_j^{(n)}$, $1\le j \le n$, be i.i.d.\ random variables such that $X_1^{(n)}$ has the same distribution as $S_1$ conditioned  on $\{ \nu_0 \le \ee^{n^{1/4}}\}$. Let $S^{(n)}_0:=0$ and $S^{(n)}_j:= X^{(n)}_1+...+X^{(n)}_j$ for $1\le j \le n$. Then
\begin{eqnarray*}
 &&\P\{ S_i \in I_{i,n}\, , \;  \nu_{i-1} \le \ee^{n^{1/4}}, \;
    \forall \, 1\le i\le n\}
    \\
 &=&[\, \P( \nu_0 \le \ee^{n^{1/4}})]^n \,
     \P \Big\{ \max_{0\le k\le n} S_k^{(n)} \le an^{1/3}, \;
     S^{(n)}_i \ge (a-\varepsilon)n^{1/3}-b(n-i)^{1/3}, \;
     \forall 1\le i \le n \Big\} .
\end{eqnarray*}

\noindent The second probability expression on the right-hand side is, according to Proposition \ref{p:mogulskii-array} (we easily check that condition (\ref{hyp-array}) is satisfied), $=\exp\{ - (1+o(1)) {\pi^2 \sigma^2\over 2} n^{1/3} \int_0^1 {\d t \over (\varepsilon+ b(1-t)^{1/3})^2}\}$, which is bounded by $\exp\{ - ({3\pi^2 \sigma^2\over 2b^2} - c_1(\varepsilon)) n^{1/3}\}$ for all sufficiently large $n$, with $c_1(\varepsilon)$ denoting a constant such that $\lim_{\varepsilon\to 0} c_1(\varepsilon) =0$. On the other hand, by (\ref{nu}), we have, for any $\eta>0$, $\E[(\nu_0)^\eta]= \E(N^\eta \sum_{i=1}^N A_i)$, which is finite (by H\"{o}lder's inequality and (\ref{H1})) if $\eta>0$ is chosen to be sufficiently small; thus $[\P( \nu_0 \le \ee^{n^{1/4}})]^n \to 1$ as $n\to \infty$. Accordingly, for all sufficiently large $n$ and some constant $c_2(\varepsilon)$ satisfying $\lim_{\varepsilon\to 0} c_2(\varepsilon) =0$,
\begin{equation}
    \P \Big( \min_{|x|=n} \overline{V}(x) \le a n^{1/3} \Big)
    \ge
    {\exp\{ n^{1/3}
    [ a - { 3 \pi^2 \sigma^2 \over 2 b^2} - c_2(\varepsilon)] \}
    \over 1+ \ee^{n^{1/4}}\sum_{j=1}^n h_{j,n}}.
    \label{ezn}
\end{equation}

We now estimate $\sum_{j=1}^n h_{j,n}$. By definition,
\begin{eqnarray*}
    h_{j, n}
 &=&\sup_{0\le u\le \varepsilon  n^{1/3} + b(n-j)^{1/3} }
    \E \Big( \ee^{ S_{n-j} }
    {\bf 1}_{ \{ u \ge  S_i \ge u - \varepsilon  n^{1/3}
    - b (n-j-i)^{1/3}, \; \forall \, i\le n-j\} }\Big)
    \\
 &\le& \sup_{0\le u\le \varepsilon  n^{1/3} + b(n-j)^{1/3} }
    \ee^u \, \P \Big( u \ge  S_i \ge u - \varepsilon  n^{1/3}
    - b (n-j-i)^{1/3}, \; \forall \, i\le n-j \Big) .
\end{eqnarray*}

\noindent Let $A$ be an integer such that $A\ge {1\over \varepsilon^2}$. Let $n_\ell:= \ell \lfloor {n \over A} \rfloor$ for $\ell =0$, $1$, $\cdots$, $A-1$ and $n_A:=n$. If $j\in [n_\ell, \, n_{\ell+1}] \cap \z$ (for some $0\le \ell \le A-1$), then
$$
h_{j, n} \le \ee^{ \varepsilon n^{1/3} + b(n-n_\ell)^{1/3} } \sup_{0\le u\le (b+\varepsilon)n^{1/3}}
    \P\Big( u\ge S_i \ge u -\varepsilon n^{1/3} -b(n-n_\ell-i)^{1/3},
    \; \forall i\le n- n_{\ell+1}\Big).
$$

\noindent We now bound the supremum on the right-hand side. If $\ell$ is such that $1- {\ell +1 \over A} \le \varepsilon$, then we simply say that the supremum is bounded by $1$, so that $\max_{n_\ell \le j \le n_{\ell+1}} h_{j, n} \le \ee^{  \varepsilon  n^{1/3} + b(n-n_\ell)^{1/3}}$. If $1- {\ell +1 \over A} > \varepsilon$, we bound the supremum by applying Corollary \ref{c:tech1} (i) to $f(t) := {\varepsilon\over (1- {\ell +1 \over A})^{1/3}} + b ( {A-\ell\over A-(\ell+1)} -t)^{1/3}$: since $f(t) \le \varepsilon^{2/3} + b(1+ {1\over \varepsilon A}-t)^{1/3}\le \varepsilon^{2/3} + b(1+ \varepsilon -t)^{1/3}$ (using $A\ge {1\over \varepsilon^2}$ for the second inequality), we have $\int_0^1 {\!\d t \over f^2(t)} \ge {3\over b^2}-c_3(\varepsilon)$, with $c_3(\varepsilon)$ denoting a constant satisfying $\lim_{\varepsilon\to 0} c_3(\varepsilon)=0$; hence by Corollary \ref{c:tech1} (i),
$$
\max_{n_\ell \le j \le n_{\ell+1}} h_{j, n} \le \ee^{\varepsilon  n^{1/3} + b(n-n_\ell)^{1/3} - ({3 \pi^2 \sigma^2 \over 2 b^2} -c_3(\varepsilon)) (n-n_{\ell+1})^{1/3}}.
$$

\noindent Therefore, for all sufficiently large $n$ and a constant $c(\varepsilon)$ satisfying $\lim_{\varepsilon\to 0} c(\varepsilon)=0$, we have, uniformly in all $\ell \in [0, \, A-1]\cap \z$,
$$
\max_{n_\ell \le j \le n_{\ell+1}} h_{j, n}\le  \ee^{n^{1/3}[( b -{ 3 \pi^2  \sigma^2 \over 2 b^2})^+ + c(\varepsilon)]},
$$

\noindent where $u^+ := \max\{ u, \, 0\}$. As a consequence, $\max_{0 \le j \le n } h_{j, n} = \max_{0\le \ell \le A-1} \max_{n_\ell \le j \le n_{\ell+1}} h_{j, n} \le \ee^{n^{1/3}[( b -{ 3 \pi^2  \sigma^2 \over 2 b^2})^+ + c(\varepsilon)]}$ for all sufficiently large $n$. In view of (\ref{ezn}), we obtain that, for any $b> a>0$,
\begin{equation}
    \liminf_{n\to \infty} \,
    {1\over n^{1/3}} \log \P \Big( \min_{|x|=n} \overline{V}(x) \le a n^{1/3} \Big)
    \ge
    -\big( b -  { 3 \pi^2  \sigma^2 \over 2 b^2} \big)^+
    + a-{3\pi^2\sigma^2 \over 2 b^2}.
    \label{V>}
\end{equation}

We now fix $a> ( { 3 \pi^2 \sigma^2 \over 2} )^{1/3}$ and $\eta>0$. We can choose $b>a$ sufficiently close to $a$ such that $( b -  { 3 \pi^2  \sigma^2 \over 2 b^2} )^+ - a+  { 3 \pi^2  \sigma^2 \over 2 b^2} < \eta$; accordingly, for all sufficiently large $n$,
\begin{equation}
    \P \Big( \min_{|x|=n} \overline{V}(x) \le an^{1/3} \Big)
    \ge
    \ee^{-\eta\, n^{1/3}}.
    \label{mcdiarmid}
\end{equation}

{F}rom here, it is routine (McDiarmid~\cite{mcdiarmid}) to obtain the upper bound in (\ref{max-ray}); we produce the details for the sake of completeness. Let $R_n := \inf\{ k: \; \# \{ x: \, |x|=k\} \ge \ee^{2\eta \, n^{1/3}}\}$. For all large $n$,
\begin{eqnarray*}
 &&\P\Big\{ R_n <\infty, \;
    \max_{k\in [{n\over 2}, \, n]}
    \min_{|x|= k+R_n} \overline{V}(x)
    > \max_{|y|=R_n} \overline{V} (y) +
    a n^{1/3} \Big\}
    \\
 &\le&\sum_{k\in [{n\over 2}, \, n]}
    \P\Big\{ R_n <\infty, \;
    \min_{|x|= k+R_n} \overline{V}(x)
    > \max_{|y|=R_n} \overline{V} (y) +
    a n^{1/3} \Big\}
    \\
 &\le&\sum_{k\in [{n\over 2}, \, n]}
    \Big[ \P\Big\{ \min_{|x| =k} \overline{V}(x)
    > a n^{1/3} \Big\}
    \Big]^{\lfloor \ee^{2\eta n^{1/3}}\rfloor}  ,
\end{eqnarray*}

\noindent which, according to (\ref{mcdiarmid}), is summable in $n$. By the Borel--Cantelli lemma, $\P$-a.s.\ for all large $n$, we have either $R_n=\infty$, or $\max_{k\in [{n\over 2}, \, n]} \min_{|x|= k+R_n} \overline{V}(x) \le \max_{|y|=R_n} \overline{V} (y) +a n^{1/3}$.

By the law of large numbers for the branching random walk (Biggins~\cite{biggins}), there exists a constant $c\in (0, \, \infty)$ such that ${1\over n} \max_{|y|=n} V(y) \to c$, $\P$-almost surely upon the system's survival. In particular, upon survival, $\max_{|y|=n} V(y)\le 2c n$, $\P$-almost surely for all large $n$. Consequently, upon the system's survival, $\P$-almost surely for all large $n$, we have either $R_n=\infty$, or $\max_{k\in [{n\over 2}, \, n]} \min_{|x|= k+R_n} \overline{V}(x) \le 2c R_n +a n^{1/3}$.

Recall that the number of particles in each generation forms a supercritical Galton--Watson process. In particular, conditionally on the system's survival, ${\# \{ u: \, |u|= k\} \over (\E N)^k}$ converges a.s.\ to a (strictly) positive random variable when $k\to \infty$, which implies $R_n \sim 2\eta {n^{1/3} \over \log (\E N)}$ $\P$-a.s.\ ($n\to \infty$), and $\max_{k\in [{n\over 2}, \, n]} \min_{|x|= k+R_n} \overline{V}(x) \ge \min_{|x|= n} \overline{V}(x)$ $\P$-almost surely for all large $n$. As a consequence, upon the system's survival, we have, $\P$-almost surely for all large $n$,
$$
\min_{|x|= n} \overline{V}(x) \le {5c \eta\over \log (\E N)} \, n^{1/3} + a n^{1/3}.
$$

\noindent Since $a$ (resp.~$\eta$) can be as close to $({3 \pi^2 \sigma^2 \over 2})^{1/3}$ (resp.~$0$) as possible, this yields the upper bound in (\ref{max-ray}), and completes the proof of Theorem \ref{t:max-ray}.\hfill$\Box$

\bigskip

Our proof of Theorem \ref{t:max-ray} gives the following deviation probability of $\min_{|x|=n} \overline{V}(x)$, which may be of independent interest.

\begin{proposition}
 \label{p:dev}
 Assume $\psi(1)=\psi'(1)=0$. For any
 $0<a\le ({ 3 \pi^2 \sigma^2 \over 2})^{1/3}$, we have
\begin{equation}
    \lim_{n \to\infty}
    {1\over n^{1/3}} \log
    \P\Big( \min_{|x|=n} \overline{V}(x) \le a n^{1/3} \Big)
    =
    a- \Big( { 3 \pi^2 \sigma^2 \over 2} \Big)^{1/3} .
    \label{V=}
\end{equation}

\end{proposition}

\noindent {\it Proof.} If $0<a<({3\pi^2 \sigma^2 \over 2})^{1/3}$, the upper and lower bounds in (\ref{V=}) follow from (\ref{V<}) and (\ref{V>}), respectively. [In (\ref{V>}), we use the fact that $b:= ({ 3 \pi^2 \sigma^2 \over 2})^{1/3}$ solves $b= {3\pi^2 \sigma^2\over 2b^2}$.]

If $a= ({ 3 \pi^2 \sigma^2 \over 2})^{1/3}$, only the lower bound in (\ref{V=}) requires a proof, which follows immediately from (\ref{V>}).\hfill$\Box$

\begin{remark}
 \label{r:zeitouni-jaffuel}
 {\rm
 Assume $\psi(1)=\psi'(1)=0$. Theorem \ref{t:max-ray} says that, on the
 set of non-extinction,
 $\P$-almost surely for $n\to \infty$,
 there exists $x_n$ with $|x_n|=n$ such that
 $\overline{V}(x_n) =
 (1+o(1)) ({ 3 \pi^2 \sigma^2 \over 2})^{1/3} n^{1/3}$. One may wonder
 whether the vertices $(x_n)$ can be chosen to form an infinite ray
 (i.e., each $x_n$ is a child of $x_{n-1}$). The answer is no:
 Jaffuel~\cite{jaffuel} proves that this is possible only
 if we increase the function
 $({ 3 \pi^2 \sigma^2 \over 2})^{1/3} n^{1/3}$
 to $({81\pi^2 \sigma^2 \over 8})^{1/3} n^{1/3}$.\hfill$\Box$
 } 
\end{remark}

\section{An estimate for one-dimensional random walks}
\label{s:spine}

$\phantom{aob}$We present in this section a probability estimate for one-dimensional random walks. It will be used in the proofs of Theorems \ref{t:main} and \ref{t:main2} in the forthcoming sections. For each $n\ge 1$, let $X_i^{(n)}$, $1\le i \le n$, be i.i.d.\ real-valued variables; let $S^{(n)}_0 :=0$ and $S^{(n)}_j := \sum_{i=1}^j X^{(n)}_i$ for $1\le j\le n$. Let $(a_n)$ be positive numbers such that $a_n \to\infty$ and ${a^2_n\over n} \to 0$, $n\to \infty$. We write $\overline{S}^{(n)}_j := \max_{1\le i\le j} S^{(n)}_i$ for $1\le j\le n$.

\begin{proposition}
 \label{p:sbar}
  Assume (\ref{hyp-array}).  Let $f: [0, 1] \to (0, \infty)$  be a continuous function.  For $\delta\ge 0$, we consider the event
  $$
  G_\delta(n)
  :=
  \Big\{  (1+\delta) \overline{S}^{(n)}_j - S^{(n)}_j \le
  a_n\, f ({j\over n}), \;
  \forall 1\le j \le n\Big\} .
  $$

 {\rm (i)} If $\delta=0$, then
 $$
 \lim_{n\to\infty} {a_n^2 \over n} \, \log  \P \Big\{  G_0(n)\Big\}  =  - {\pi^2 \sigma^2\over 8} \int_0^1 {\d s \over f^2(s)} .
 $$
 Moreover, for any fixed  $0<b<1$,
 $$
 \lim_{n\to\infty} {a_n^2 \over n} \, \log  \P \Big\{
  G_0(n), \, \overline{S}^{(n)}_n - S^{(n)}_n \le b \, a_n\, f(1)  \Big\}  =  - {\pi^2 \sigma^2\over 8} \int_0^1 {\d s \over f^2(s)} .
  $$

 {\rm (ii)} If $\delta>0$, then
 $$
 \lim_{n\to\infty} {a_n^2 \over n} \, \log  \P \Big\{  G_\delta(n)\Big\}  =  -  {\pi^2 \sigma^2\over 2 } \int_0^1 {\d s \over f^2(s)} .
 $$
\end{proposition}

We mention that for the centered random walk $(S_n)$ given in (\ref{change-proba}), assumption (\ref{hyp-array}) is obviously satisfied. Hence Proposition \ref{p:mogulskii-array} as well as Corollary \ref{c:tech} below, hold also for $(S_n)$.

\bigskip

\noindent {\it Proof of Proposition \ref{p:sbar}.} (i) Let $0<\varepsilon< {1\over 4} \min \{ b, \, \min_{0\le t \le 1} f(t)\}$ and let $A$ be a large integer. Consider a sufficiently large $n$ such that $ \sup_{0\le s< t\le 1, \, t-s \le 2 A a_n^2/n} |f(t) - f(s)|\le \varepsilon$. Let $m= \lfloor {n \over A^2 a_n^2}\rfloor $. For $0 \le k < A m$, let $r_k:= k \lfloor A a_n^2 \rfloor$ and $r_{A m}:= n$. Note that $\lfloor A a_n^2 \rfloor \le r_{Am } - r_{A m-1} \le 2 \lfloor A a_n^2 \rfloor$. Let $\ell \in [0, \, A-1]\cap \z$ and $k \in [ \ell m, \, (\ell +1)m -1] \cap \z$. For all $r_k \le j < r_{k+1}$, $|f({j \over n}) - f( {\ell \over A})|\le  \varepsilon $. Define
$$
E_n^{(\pm)}
:=
\bigcap_{\ell=0}^{A-1} \,
\bigcap_{k= \ell m}^{ (\ell+1) m - 1} \,
\bigcap_{j=r_k}^{r_{k+1}-1}
\Big\{ \overline{S}^{(n)}_{j} - S^{(n)}_{j } \le a_n\,  (f ( {\ell
\over A} ) \pm \varepsilon )\Big\}.
$$

\noindent Then
\begin{eqnarray*}
    \P\Big( G_0(n)\Big)
 &\le &\P\left(E_n^{(+)}\right) ,
    \\
    \P\Big( G_0(n),  \overline{S}^{(n)}_n - S^{(n)}_n
    \le b \, a_n\, f ({j\over n}) \Big)
 &\ge & \P\Big( E^{(-)}_n \cap
    \bigcap_{0\le k \le A m}
    \{ \overline{S}^{(n)}_{r_k} - S^{(n)}_{r_k} \le \varepsilon a_n \}
    \Big) .
\end{eqnarray*}

\noindent Observe that for any $r_k$, conditionally on
$\sigma\{S^{(n)}_j, 0\le j\le r_k \}$ and on $\{ \overline
S^{(n)}_{r_k} - S^{(n)}_{r_k}= x\}$, the reflecting process $(
\overline{S}^{(n)}_{i+r_k} -S^{(n)}_{i+r_k}, \; 0\le i \le r_{k+1}-
r_k)$ has the same law as $(\max\{ x, \, \overline{S}^{(n)}_i\} -
S^{(n)}_i, \; 0\le i \le r_{k+1}- r_k)$. Using this observation for all $k$, we see that
\begin{eqnarray}
 && \P(E_n^{(+)}) \le
    \prod_{\ell=0}^{A-1} \prod_{k= \ell m}^{ (\ell+1) m -1} \,
    \P\Big\{ \max_{0\le i<r_{k+1}- r_k}
    (\overline{S}^{(n)}_i - S^{(n)}_i)\le
    a_n\, (f( {\ell \over A} ) +\varepsilon ) \Big\} ,
    \label{pen+}
    \\
 && \P\Big( E^{(-)}_n \cap
    \bigcap_{0\le k \le A m}
    \{ \overline{S}^{(n)}_{r_k} - S^{(n)}_{r_k}
    \le  \varepsilon a_n \}\Big)
    \ge \prod_{\ell=0}^{A-1}
    \prod_{k= \ell m}^{ (\ell+1) m - 1}
    \, \P\Big\{   \Upsilon_k \Big\} ,
   \label{pen-}
\end{eqnarray}

\noindent with
$$
\Upsilon_k:=\Big\{ \max_{0\le i<r_{k+1}- r_k}(
\overline{S}^{(n)}_i - S^{(n)}_i ) \le a_n\,  (f ( {\ell \over A} ) -2
\varepsilon ), \overline{S}^{(n)}_{r_{k+1}- r_k} - S^{(n)}_{r_{k+1}-
r_k} < \varepsilon a_n, \, \overline{S}^{(n)}_{r_{k+1}- r_k} >
\varepsilon a_n\Big\}.
$$

Now, we prove the upper bound in (i). By (\ref{pen+}),
$$
{a_n^2\over n } \log \P(E_n^{(+)}) \le {m a_n ^2 \over n
}\sum_{\ell=0}^{A-1} \log \P\Big\{ \overline{S}^{(n)}_i - S^{(n)}_i
\le a_n \,  (f ( {\ell \over A} ) +\varepsilon ), \,\forall\, 0\le
i<\lfloor A a_n^2 \rfloor\Big\}.
$$

\noindent According to Donsker's invariance principle,\footnote{Finite-dimensional convergence is checked by Lindeberg's condition in the central limit theorem, whereas tightness is proved via a standard argument as in Billingsley~\cite{billingsley}.} the probability term on the right-hand side converges, when $n\to \infty$, to
$$
\P  \Big\{ \sup_{0\le t\le 1}  ( \overline{W}(t) - W(t) ) \le {1\over \sigma \sqrt{A}}  (f ( {\ell \over A}) +\varepsilon )\Big\},
$$

\noindent where $W$ is a standard one-dimensional
Brownian motion, and $\overline{W}(t)= \sup_{0\le
s\le t} W(s)$. By L\'evy's identity, $(\overline{W}(t) - W(t), \; t\ge 0)$ is distributed as $(|W(t)|, \; t\ge 0)$; thus we have
\begin{equation}
    \P\Big\{ \sup_{0\le t\le 1} ( \overline{W}(t) - W(t) ) \le u\Big\}
    =
    \ee^{ - (1+o(1)) {\pi^2 \over 8u^2}}, \qquad u\to 0,
    \label{chung}
\end{equation}

\noindent which can be easily deduced from Formula (5.9) of page 342 of Feller~\cite{Feller}, taking $a=2u$, $t=1$ and $x=u$. As a consequence, for all sufficiently large $A$, say $A \ge A_0= A_0(\varepsilon, \, \sigma, \, f)$,
$$
\log \P\Big\{ \sup_{0\le t\le 1} ( \overline{W}(t) - W(t) ) \le {1\over \sigma \sqrt{A}}  (f ( {\ell \over A} ) +\varepsilon )\Big\}
\le
- {(1-\varepsilon)\pi^2 \sigma^2 A \over 8 (f ( {\ell \over A}) +\varepsilon )^2}.
$$

\noindent  Since $m \sim  {n\over a_n^2 A^2}$, we get, for $A \ge A_0$,
\begin{eqnarray*}
    \limsup_{n \to\infty}{a_n^2\over n } \log \P(E_n^{(+)})
 &\le& {1\over A^2} \sum_{\ell=0}^{A-1}
    \log \P\Big\{ \sup_{0\le t\le1}  ( \overline{W}(t) - W(t) ) \le
    {1\over \sigma \sqrt{A}} (f ( {\ell \over A} ) +\varepsilon )\Big\}
    \\
 &\le& - {\pi^2 \sigma^2 \over 8 } {1- \varepsilon\over A}
    \sum_{\ell=0}^{A-1}
    {1\over (f ( {\ell \over A}) +\varepsilon )^2} .
\end{eqnarray*}

\noindent Letting $A\to\infty$ and then $\varepsilon\to0$, we get
the upper  bound in (i):
$$
\limsup_{n\to\infty} {a_n^2 \over n} \, \log  \P \Big\{ \overline
S^{(n)}_i - S^{(n)}_j \le a_n\, f({j\over n}), \,\, \forall 1\le j
\le n\Big\} \le  - \, {\pi^2 \sigma^2\over 8} \, \int_0^1 {\d s \over
f^2(s)} .
$$

To prove the lower bound in (i), we go back to the events
$\Upsilon_k$ in (\ref{pen-}). Observe that for each
$1\le i \le r_{k+1}- r_k$, all the three events in $\Upsilon_k$ are
non-decreasing with respect to $S^{(n)}_i-S^{(n)}_{i-1}$. By the FKG
inequality,
\begin{eqnarray*}
    \P\Big( \Upsilon_k \Big)
  &\ge& \P\Big(
    \max_{0\le i<r_{k+1}- r_k}( \overline{S}^{(n)}_i
    - S^{(n)}_i ) \le a_n\, (f( {\ell \over A}) -2 \varepsilon)\Big)
   \\
 &&\qquad \times
    \P\Big( \overline{S}^{(n)}_{r_{k+1}- r_k}
    - S^{(n)}_{r_{k+1}- r_k} < \varepsilon a_n\Big) \,
   \P\Big( \overline{S}^{(n)}_{r_{k+1}- r_k} > \varepsilon a_n\Big).
\end{eqnarray*}

\noindent Recall that $r_{k+1}-r_k = \lfloor A a_n^2\rfloor$ for $0\le k< Am-1$, and $\lfloor A a_n^2\rfloor \le r_{Am} -r_{Am-1} \le 2 \lfloor A a_n^2\rfloor$. Using Donsker's invariance principle again, we see that there exists a constant $c(\varepsilon)>0$ such that for all $k$, $\P (\overline{S}^{(n)}_{r_{k+1}- r_k} - S^{(n)}_{r_{k+1}- r_k} <
\varepsilon a_n ) \, \P ( \overline{S}^{(n)}_{r_{k+1}- r_k} >
\varepsilon a_n ) \ge c(\varepsilon)$. {F}rom this, the lower
bound in (i) follows in the same way as the upper bound in (i).

(ii) Let us first prove the following fact: for any  fixed $c
>0$, 
\begin{equation}
    \label{toro1} 
    \P\Big\{ \sup_{0\le s\le 1} (\overline{W}(s) - W(s))
\le
    u,  \overline W(1) \le c \, u\Big\}
    =
    \ee^{  - {\pi^2\over 2 u^2} (1+o(1))}, \qquad u\to 0.
\end{equation}

To see why (\ref{toro1}) holds, we denote by $L(t)$ the local time
at $0$ of $W$ up to time $t$, and recall from Borodin and Salminen
(\cite{BS02}, page 259, Formula 1.16.2) that, for $\lambda>0$,
$$
\int_0^\infty \ee^{-\lambda t}\, \P\Big(\sup_{s\le t}   |W(s)| \le  1, \, L(t) \le c \Big) \d t = {1\over \lambda} \, \Big( 1- {1\over \cosh( \sqrt{2\lambda})}\Big) \, \Big( 1- \ee^{ - c \sqrt{\lambda\over 2}\, \coth( \sqrt{2\lambda})}\Big).
$$

\noindent By analytic continuation, we get that for $0< \lambda<
{\pi^2\over 2}$,
$$
\int_0^\infty \ee^{\lambda t} \, \P\Big(\sup_{s\le t}   |W(s)| \le  1, \, L(t) \le c \Big) \d t = {1\over \lambda} \, \Big(  {1\over \cos( \sqrt{2\lambda})}-1\Big) \, \Big( 1- \ee^{ - c\,  \sqrt{\lambda\over 2}\, \mbox{cotan}( \sqrt{2\lambda})}\Big).
$$


\noindent This implies, by means of a Tauberian theorem (see, for
example, Theorem 3.2 of \cite{HS97}), that
$$
\P\Big(\sup_{0\le s\le t} |W(s)| \le 1, \, L(t) \le c\Big) =\ee^{ - ({\pi^2\over 2} +o(1)) t}, \qquad t\to\infty ,
$$

\noindent which, by scaling, is equivalent to $\P (\sup_{0\le s\le
1} |W(s)| \le u, \, L(1) \le {u\over \delta}) =\ee^{ - ({\pi^2\over
2u} +o(1)) }$, $u\to 0$.

\noindent By L\'evy's identity, the two processes $(\overline{W} -
W, \, \overline{W})$ and $(|W|, \, L)$ have the same law;
consequently, this implies (\ref{toro1}).

Now let us proceed to prove the upper bound in (ii). Let $\varepsilon>0$, and let $(r_k)$ be as in the proof of (i), i.e., $A$ is a large integer, $m:= \lfloor {n \over A^2 a_n^2}\rfloor $, $r_k:= k \lfloor A a_n^2 \rfloor$ (for $0 \le k < A m$) and $r_{A m}:= n$, with $n$ sufficiently large such that $|f({j \over n}) - f( {\ell \over A})|\le \varepsilon$ for $r_k \le j \le r_{k+1}$ and $k\in [\ell m, \, (\ell +1)m) \cap \z$. Let
$$ 
F_n^{(+)} 
:=
\bigcap_{\ell=0}^{A-1} \, \bigcap_{k= \ell m}^{ (\ell+1) m - 1} \,
\Big\{ \max_{ r_k\le j \le r_{k+1}} (\overline{S}^{(n)}_{j} -
S^{(n)}_{j }) \le a_n\, (f ( {\ell \over A} ) +\varepsilon) , \; 
\overline{S}^{(n)}_{r_{k+1}} \le  {c\over \delta}\, a_n\Big\},
$$

\noindent where $c=\sup_{0\le t \le 1} f(t)$. Clearly $G_\delta(n)
\subset F_n^{(+)}$. The Markov property yields that for each $k$,
conditionally on $\sigma\{S^{(n)}_j, \; 0\le j\le r_k \}$ and on $\{
\overline S^{(n)}_{r_k} - S^{(n)}_{r_k}= x_k, \; S^{(n)}_{r_k}=y_k\}$, the process $( \overline{S}^{(n)}_{i+r_k}, \, \overline{S}^{(n)}_{i+r_k} -S^{(n)}_{i+r_k}, \; 0\le i \le r_{k+1}- r_k)$ has the same law as $(\max\{x_k, \, \overline{S}^{(n)}_i\}+y_k,\,  \max\{ x_k, \,
\overline{S}^{(n)}_i\} - S^{(n)}_i, \; 0\le i \le r_{k+1}- r_k)$.
On $F_n^{(+)}$, we have $-y_k \le x_k \le (c+\varepsilon) a_n$ (recalling that $c= \sup_{0\le t \le 1} f(t)$), thus $y_k \ge - (c+\varepsilon) a_n$. Therefore, by the Markov property,
$$ 
\P(F_n^{(+)}) 
\le
    \prod_{\ell=0}^{A-1} \prod_{k= \ell m}^{ (\ell+1) m -1} \,
    \P\Big\{ \max_{0\le i\le r_{k+1}- r_k}
    (\overline{S}^{(n)}_i - S^{(n)}_i)\le
    a_n\, (f( {\ell \over A} )+\varepsilon)  , \, \overline{S}^{(n)}_{r_{k+1}-r_k} \le   (c+\varepsilon+ {c \over \delta}) a_n \Big\} .
$$

This is the analogue of (\ref{pen+}) for (ii). From here, the rest of the proof of the upper bound in (ii) is done by using exactly the same arguments as in (i), by applying (\ref{toro1}) instead of (\ref{chung}). We omit the details.

The proof of the lower bound in (ii) is easy. Indeed, let $0<\varepsilon< \inf_{t\in [0, \, 1]}f(t)$, and let
$$ 
F^{(-)}_n
:= 
\left\{ - ( f ({i\over n}) - \varepsilon) \le
{S^{(n)}_i\over a_n} \le { \varepsilon\over 1+\delta}, \, \forall 0
\le i \le n\right\}.$$

\noindent Clearly $F^{(-)}_n \subset G_\delta(n)$. By (\ref{mogulskii-array}), we have 
$$ 
\lim_{n \to\infty} {a_n^2 \over n} \log \P ( F^{(-)}_n ) 
= 
- {\pi^2 \sigma^2 \over 2} \int_0^1 {\d t \over ( f(t) - \varepsilon +
{\varepsilon\over 1+\delta})^2}.$$

\noindent Letting $ \varepsilon\to 0$ gives the lower bound in (ii).\hfill$\Box$

\bigskip

The following corollary follows from Proposition \ref{p:sbar} exactly as Corollary \ref{c:tech1} follows from Proposition \ref{p:mogulskii-array}.

\begin{corollary}
 \label{c:tech}
 Assume that (\ref{hyp-array}) is satisfied with $a_n = n^{1/3}$. Let
 $a>0$ and $\delta>0$. Then for $n\to \infty$,
 \begin{eqnarray*}
     \sum_{j=1}^n  \ee^{-a (n-j)^{1/3}}
     \P\Big( \overline{S}^{(n)}_i -  S^{(n)}_i \le
     a (n-i)^{1/ 3}, \, \forall\, 1\le i \le j\Big)
  &=& \ee^{- \min\{ a, \, {3\pi^2 \sigma^2 \over 8 a^2}\}
     (1+o(1)) n^{1/3}} ,
     \\
     \sum_{j=1}^n  \ee^{-a (n-j)^{1/3}}
     \P\Big((1+\delta) \overline{S}^{(n)}_i
     -  S^{(n)}_i \le a (n-i)^{1/ 3}, \, \forall\, 1\le i \le j\Big)
  &=& \ee^{- \min\{ a, \,
     {3\pi^2 \sigma^2 \over 2  a^2}\} (1+o(1)) n^{1/3}}.
 \end{eqnarray*}

\end{corollary}

\section{Proof of Theorem \ref{t:main2}}
\label{s:proof-psi'(1)>0}

$\phantom{aob}$We assume $\inf_{t\in [0, \, 1]} \psi(t)=0$ and $\psi'(1) \ge 0$ in this section. Let $\theta\in (0, \, 1]$ be such that $\psi'(\theta)=0$ as in (\ref{theta}). By Theorem \ref{t:max-ray} and
(\ref{hitting-proba<=>potential}), we get that, on the set of non-extinction,
$$
\liminf_{n\to \infty} \, {1\over n^{1/3}} \log \varrho_n
\ge
- \alpha_\theta^{1/3},
\qquad\hbox{\rm $\P$-a.s.},
$$

\noindent where $\alpha_\theta := {3\pi^2\over 2\theta} \E[ \sum_{i=1}^N A_i^\theta (\log A_i)^2] = {3\pi^2\over 2\theta} \E[ \sum_{|x|=1} V(x)^2 \ee^{-\theta V(x)}]$, and $\varrho_n := P_\omega \{ \tau_n < \tau_0\}$ is as in (\ref{hitting-proba<=>potential}). In view of Fact \ref{f:rhonxk}, it remains only to check that if $\psi'(1)>0$ (i.e., if $\theta<1$), then we have, on the set of non-extinction,
\begin{equation}
    \limsup_{n\to \infty} \, {1\over n^{1/3}} \log \varrho_n
    \le
    - \alpha_\theta^{1/3},
    \qquad\hbox{\rm $\P$-a.s.}
    \label{upperrhon}
\end{equation}

We do not assume $\psi'(1)>0$ for the moment (so $\theta$ can be $1$, and the inequality (\ref{rho-ub}) below can also be used in the proof of Theorem \ref{t:main} in the next section). Let $a>0$, $n\ge 1$ and $\delta \ge 0$. For any $y$ with $|y|\le n$, say $|y|=j$, we introduce the following event:
$$
E_\delta(y) = \Big\{  (1+\delta) \overline{V}(y) - V(y) \ge {a\over
\theta} (n-j)^{1/3}\Big\} \cap \bigcap_{i=1}^{j-1} \Big\{
(1+\delta)\overline{V}(y_i) - V(y_i) < {a\over \theta}
(n-i)^{1/3} \Big\},
$$

\noindent where $y_i$ is the unique vertex of $[\![ \varnothing, \, y]\!]$ that is in the $i$-th generation, whereas $\overline{V}(x) := \max_{z\in \, ]\!]\varnothing, \, x]\!]} V(z)$ as in (\ref{Vbar}).

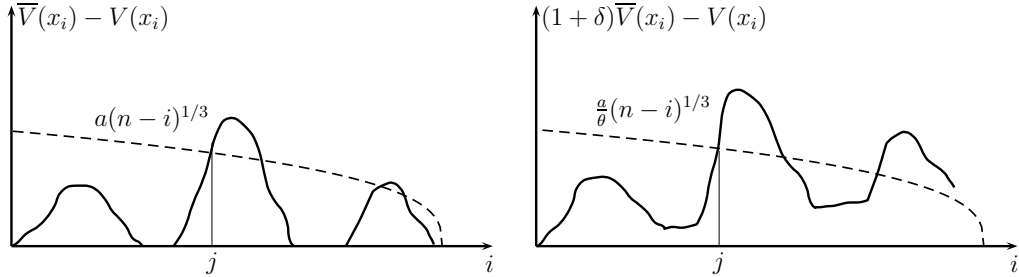
\begin{figure}[h]

\centering

\scalebox{0.8} 
{
\begin{pspicture}(0,-2.3129)(17.78,2.3189)
\rput(0.84,-1.6871){\psaxes[linewidth=0.04,labels=none,ticks=none,ticksize=0.025799999cm]{->}(0,0)(0,0)(8,4)}
\pscustom[linewidth=0.04]
{
\newpath
\moveto(0.86,-1.6811)
\lineto(0.89726496,-1.6602663)
\curveto(0.91589737,-1.6498494)(0.95316243,-1.6185994)(0.9717949,-1.5977663)
\curveto(0.99042743,-1.5769331)(1.0230342,-1.540475)(1.0370085,-1.52485)
\curveto(1.0509828,-1.509225)(1.0882479,-1.477975)(1.1115385,-1.46235)
\curveto(1.1348292,-1.446725)(1.1907265,-1.415475)(1.2233334,-1.39985)
\curveto(1.2559401,-1.384225)(1.33047,-1.321725)(1.3723933,-1.27485)
\curveto(1.4143164,-1.227975)(1.4655555,-1.165475)(1.4748716,-1.14985)
\curveto(1.4841878,-1.134225)(1.5074785,-1.0977665)(1.5214529,-1.0769331)
\curveto(1.5354275,-1.0560997)(1.5587181,-1.0144331)(1.5680342,-0.9936)
\curveto(1.5773504,-0.9727669)(1.5913247,-0.9311003)(1.5959828,-0.9102669)
\curveto(1.6006413,-0.88943344)(1.6192737,-0.8477666)(1.6332479,-0.82693315)
\curveto(1.647222,-0.8060997)(1.684487,-0.76443315)(1.7077779,-0.7436)
\curveto(1.7310686,-0.7227669)(1.7776498,-0.6967253)(1.8009402,-0.6915169)
\curveto(1.8242306,-0.68630844)(1.8894441,-0.6811)(1.9313675,-0.6811)
\curveto(1.9732909,-0.6811)(2.0431626,-0.6811)(2.071111,-0.6811)
\curveto(2.0990596,-0.6811)(2.1596153,-0.696725)(2.1922224,-0.71235)
\curveto(2.2248294,-0.727975)(2.280727,-0.7696416)(2.304017,-0.79568315)
\curveto(2.3273072,-0.8217247)(2.3738885,-0.8790163)(2.3971794,-0.9102666)
\curveto(2.4204705,-0.9415169)(2.457735,-1.0040169)(2.4717093,-1.0352665)
\curveto(2.4856832,-1.0665163)(2.5229483,-1.1185997)(2.5462391,-1.1394331)
\curveto(2.5695305,-1.1602665)(2.625428,-1.2175584)(2.6580343,-1.2540169)
\curveto(2.6906404,-1.2904752)(2.737222,-1.3477669)(2.7511964,-1.3686)
\curveto(2.7651713,-1.3894331)(2.7884622,-1.4363081)(2.7977777,-1.46235)
\curveto(2.8070931,-1.4883919)(2.835042,-1.5352669)(2.8536751,-1.5561)
\curveto(2.8723085,-1.5769331)(2.9142318,-1.602975)(2.9375215,-1.6081837)
\curveto(2.9608114,-1.6133925)(2.9980764,-1.6528505)(3.0120513,-1.6870998)
}
\pscustom[linewidth=0.04]
{
\newpath
\moveto(6.4,-1.6811)
\lineto(6.48,-1.532615)
\curveto(6.52,-1.4583724)(6.5800004,-1.3470088)(6.6000004,-1.3098879)
\curveto(6.6200004,-1.2727668)(6.685,-1.1348882)(6.73,-1.0341305)
\curveto(6.775,-0.9333726)(6.8300004,-0.81140274)(6.84,-0.7901909)
\curveto(6.8500004,-0.768979)(6.915,-0.72125167)(6.9700003,-0.6947365)
\curveto(7.025,-0.6682213)(7.105,-0.636403)(7.13,-0.6311)
\curveto(7.155,-0.625797)(7.21,-0.6417061)(7.2400002,-0.6629182)
\curveto(7.27,-0.6841304)(7.32,-0.7424638)(7.34,-0.7795851)
\curveto(7.36,-0.81670636)(7.395,-0.9121608)(7.4100003,-0.9704939)
\curveto(7.425,-1.0288271)(7.46,-1.1136754)(7.48,-1.1401906)
\curveto(7.5,-1.1667061)(7.5750003,-1.2409487)(7.63,-1.2886757)
\curveto(7.685,-1.3364028)(7.75,-1.4159483)(7.76,-1.4477665)
\curveto(7.77,-1.4795848)(7.795,-1.532615)(7.81,-1.553827)
\curveto(7.8250003,-1.5750393)(7.8500004,-1.6280696)(7.86,-1.6598878)
}
\pscustom[linewidth=0.04]
{
\newpath
\moveto(3.54,-1.6811)
\lineto(3.6399999,-1.5311)
\curveto(3.6899998,-1.4561)(3.7649999,-1.2811)(3.79,-1.1811)
\curveto(3.8149998,-1.0811)(3.8899999,-0.8811)(3.9399998,-0.7811)
\curveto(3.9899998,-0.6811)(4.065,-0.4811)(4.0899997,-0.3811)
\curveto(4.115,-0.2811)(4.165,-0.0811)(4.19,0.0189)
\curveto(4.2149997,0.1189)(4.29,0.2939)(4.3399997,0.3689)
\curveto(4.39,0.4439)(4.515,0.4689)(4.5899997,0.4189)
\curveto(4.665,0.3689)(4.79,0.2439)(4.8399997,0.1689)
\curveto(4.89,0.0939)(4.965,-0.1061)(4.9900002,-0.23110001)
\curveto(5.0150003,-0.3561)(5.065,-0.5811)(5.09,-0.6811)
\curveto(5.1150002,-0.7811)(5.215,-0.9561)(5.29,-1.0311)
\curveto(5.3650002,-1.1061)(5.44,-1.3061)(5.44,-1.4311)
\curveto(5.44,-1.5561)(5.4900002,-1.6811)(5.54,-1.6811)
}
\psline[linewidth=0.01cm](4.172,-1.6871)(4.172,0.03)
\usefont{T1}{ptm}{m}{n}
\rput(3.1742,0.433){$a(n-i)^{1/3}$}
\usefont{T1}{ptm}{m}{n}
\rput(8.79,-1.9761){$i$}
\usefont{T1}{ptm}{m}{n}
\rput(4.1792,-1.9761){$j$}
\usefont{T1}{ptm}{m}{n}
\rput(2.2,2.1239){$\overline{V}(x_i)-V(x_{i})$}
\rput(9.56,-1.6871){\psaxes[linewidth=0.04,labels=none,ticks=none,ticksize=0.025799999cm]{->}(0,0)(0,0)(8,4)}
\pscustom[linewidth=0.04]
{
\newpath
\moveto(9.576624,-1.6911299)
\lineto(9.61068,-1.6653831)
\curveto(9.627708,-1.6525099)(9.6603365,-1.6164452)(9.675936,-1.593254)
\curveto(9.691536,-1.5700628)(9.718836,-1.529478)(9.7305355,-1.5120845)
\curveto(9.742235,-1.4946907)(9.774863,-1.4586263)(9.795792,-1.439955)
\curveto(9.816719,-1.4212841)(9.867805,-1.3826647)(9.897961,-1.3627162)
\curveto(9.928117,-1.3427678)(9.993372,-1.2706386)(10.028473,-1.2184578)
\curveto(10.063572,-1.166277)(10.105758,-1.0973414)(10.112843,-1.0805867)
\curveto(10.119927,-1.0638319)(10.137999,-1.0245243)(10.148985,-1.0019716)
\curveto(10.159971,-0.97941846)(10.17733,-0.9349519)(10.183701,-0.9130382)
\curveto(10.190071,-0.8911244)(10.198201,-0.8479353)(10.199959,-0.8266597)
\curveto(10.201717,-0.80538434)(10.214458,-0.76155627)(10.225445,-0.7390034)
\curveto(10.23643,-0.71645063)(10.267631,-0.6700679)(10.287844,-0.64623797)
\curveto(10.308058,-0.622408)(10.350628,-0.59022534)(10.372984,-0.5818724)
\curveto(10.39534,-0.5735194)(10.459225,-0.5594182)(10.500752,-0.5536696)
\curveto(10.54228,-0.5479213)(10.611491,-0.5383406)(10.639175,-0.5345082)
\curveto(10.66686,-0.5306758)(10.728985,-0.5378499)(10.763428,-0.5488564)
\curveto(10.79787,-0.5598628)(10.858951,-0.5934712)(10.885592,-0.61607313)
\curveto(10.912234,-0.6386751)(10.966231,-0.6890382)(10.9935875,-0.71679944)
\curveto(11.020943,-0.7445607)(11.066426,-0.8013606)(11.084553,-0.83039933)
\curveto(11.102679,-0.859438)(11.146734,-0.9059193)(11.172662,-0.92336255)
\curveto(11.198589,-0.9408055)(11.261816,-0.9898915)(11.299113,-1.0215346)
\curveto(11.3364105,-1.0531776)(11.3904085,-1.1035411)(11.407107,-1.1222612)
\curveto(11.423807,-1.1409813)(11.453304,-1.1842201)(11.466103,-1.2087386)
\curveto(11.478902,-1.2332569)(11.513014,-1.2758567)(11.534327,-1.2939382)
\curveto(11.55564,-1.3120195)(11.6007395,-1.3320671)(11.624523,-1.3340334)
\curveto(11.648308,-1.3359993)(11.690631,-1.3699751)(11.70917,-1.4019845)
}
\pscustom[linewidth=0.04]
{
\newpath
\moveto(15.064295,-0.93148607)
\lineto(15.123179,-0.773434)
\curveto(15.152621,-0.6944081)(15.196785,-0.57586926)(15.211506,-0.53635633)
\curveto(15.226227,-0.4968434)(15.271708,-0.3513545)(15.302466,-0.24537796)
\curveto(15.333226,-0.13940139)(15.370981,-0.0110420175)(15.377978,0.011340795)
\curveto(15.384975,0.033723608)(15.442817,0.08991257)(15.493662,0.12371903)
\curveto(15.544506,0.1575255)(15.619388,0.20001267)(15.643424,0.2086937)
\curveto(15.66746,0.21737473)(15.724123,0.20915726)(15.756749,0.19225907)
\curveto(15.789373,0.17536087)(15.8469,0.124434665)(15.871801,0.09040635)
\curveto(15.896703,0.056378026)(15.944461,-0.033375695)(15.967318,-0.0891011)
\curveto(15.990174,-0.1448262)(16.036478,-0.22407424)(16.059925,-0.24759658)
\curveto(16.083372,-0.2711189)(16.167843,-0.33437636)(16.228868,-0.37411088)
\curveto(16.289894,-0.41384566)(16.365187,-0.48372695)(16.379454,-0.51387346)
\curveto(16.393723,-0.54401994)(16.425758,-0.59312147)(16.443525,-0.6120763)
\curveto(16.461292,-0.63103104)(16.493328,-0.68013257)(16.507595,-0.7102791)
}
\pscustom[linewidth=0.04]
{
\newpath
\moveto(12.231309,-1.3236488)
\lineto(12.309797,-1.1614698)
\curveto(12.349041,-1.0803804)(12.399337,-0.8968806)(12.410389,-0.79447037)
\curveto(12.421441,-0.6920604)(12.468308,-0.48381484)(12.504125,-0.37797928)
\curveto(12.53994,-0.2721437)(12.586807,-0.063897856)(12.597859,0.038512427)
\curveto(12.6089115,0.14092271)(12.631014,0.34574297)(12.642066,0.44815293)
\curveto(12.653118,0.55056304)(12.703414,0.7340627)(12.742659,0.8151521)
\curveto(12.781902,0.8962416)(12.902293,0.9381157)(12.98344,0.8989)
\curveto(13.064589,0.8596843)(13.205547,0.7530812)(13.26536,0.6856938)
\curveto(13.325171,0.61830646)(13.426886,0.43061417)(13.46879,0.310309)
\curveto(13.5106945,0.19000383)(13.591072,-0.025860133)(13.629549,-0.12141921)
\curveto(13.668026,-0.2169783)(13.791078,-0.37649918)(13.875652,-0.44046095)
\curveto(13.960228,-0.5044228)(14.061943,-0.6921153)(14.0790825,-0.8158461)
\curveto(14.096224,-0.93957686)(14.162891,-1.0564564)(14.212419,-1.0496053)
}
\psline[linewidth=0.01cm](12.6,-1.7)(12.6,-0.1)
\usefont{T1}{ptm}{m}{n}
\rput(11.4942,0.573){$\frac{a}{\theta}(n-i)^{1/3}$}
\usefont{T1}{ptm}{m}{n}
\rput(17.51,-1.9761){$i$}
\usefont{T1}{ptm}{m}{n}
\rput(12.6,-1.9761){$j$}
\usefont{T1}{ptm}{m}{n}
\rput(11.52,2.1239){$(1+\delta)\overline{V}(x_i)-V(x_{i})$}
\pscustom[linewidth=0.04]
{
\newpath
\moveto(11.6742,-1.4011)
\lineto(11.7142,-1.3911)
\curveto(11.7342,-1.3861)(11.7742,-1.3761)(11.7942,-1.3711)
\curveto(11.8142,-1.3661)(11.8592,-1.3661)(11.8842,-1.3711)
\curveto(11.9092,-1.3761)(11.9542,-1.3811)(11.9742,-1.3811)
\curveto(11.9942,-1.3811)(12.0342,-1.3761)(12.0542,-1.3711)
\curveto(12.0742,-1.3661)(12.1142,-1.3561)(12.1342,-1.3511)
\curveto(12.1542,-1.3461)(12.1892,-1.3361)(12.2342,-1.3211)
}
\pscustom[linewidth=0.04]
{
\newpath
\moveto(14.1942,-1.0411)
\lineto(14.2442,-1.0411)
\curveto(14.2692,-1.0411)(14.3192,-1.0361)(14.3442,-1.0311)
\curveto(14.3692,-1.0261)(14.4192,-1.0211)(14.4442,-1.0211)
\curveto(14.4692,-1.0211)(14.5142,-1.0161)(14.5342,-1.0111)
\curveto(14.5542,-1.0061)(14.5942,-0.9961)(14.6142,-0.9911)
\curveto(14.6342,-0.9861)(14.6742,-0.9761)(14.6942,-0.9711)
\curveto(14.7142,-0.9661)(14.7592,-0.9611)(14.7842,-0.9611)
\curveto(14.8092,-0.9611)(14.8592,-0.9611)(14.8842,-0.9611)
\curveto(14.9092,-0.9611)(14.9542,-0.9561)(14.9742,-0.9511)
\curveto(14.9942,-0.9461)(15.0292,-0.9361)(15.0742,-0.9211)
}
\psplot[plotpoints=200,linestyle=dashed]{0.8658}{8}{ 8 x sub 1 3 div exp 1.701 sub}
\psplot[plotpoints=200,linestyle=dashed]{9.6658}{17}{ 17 x sub 1 3 div exp 1.701 sub}
\end{pspicture}
}
\caption{$j=\inf\{i:E_{\delta}(x_{i}) \text{ holds.}\}$}
\end{figure}

Let as before $\tau_n := \inf \{ i\ge 1: \, |X_i| =n\}$
and $T(x):= \inf \{ k\ge 0: \, X_k =x\}$. Consider any vertex $x$ with $|x|=n$. Let $j= j(x)\in [1, \, n]\cap \z$ be the smallest integer such that $(1+\delta)\overline{V}(x_j)-V(x_j) \ge {a\over \theta} (n-j)^{1/3}$. Such a $j$ exists. Moreover, we have $T(x) \ge T(x_j)$, and $E_\delta(x_j)$ holds. Consequently,
$$
\tau_n = \inf_{|x|=n} T(x)\ge \min_{1\le j\le n} \inf \{  T(y):\,
|y|=j  \mbox { and $E_\delta(y)$ holds}\} ,
$$

\noindent so that $\varrho_n =  P_\omega\{ \tau_n < \tau_0\} \le \sum_{j=1}^n \sum_{|y|=j} {\bf 1}_{E_\delta(y)} \, P_\omega \{ T(y)< \tau_0\}$. By (\ref{zk-probatrans}), we obtain:
\begin{equation}
    \varrho_n
    \le
    \sum_{j=1}^n \sum_{|y|=j} {\bf 1}_{E_\delta(y)} \,
    \omega(\varnothing, \, y_1) \, \ee^{V(y_1) - \overline{V}(y)}
    =
    \omega(\varnothing, \, {\buildrel \leftarrow \over \varnothing})
    \sum_{j=1}^n \sum_{|y|=j} {\bf 1}_{E_\delta(y)} \,
    \ee^{- \overline{V}(y)} ,
    \label{rho-ub}
\end{equation}

\noindent which is bounded by $\sum_{j=1}^n \sum_{|y|=j} {\bf 1}_{E_\delta(y)} \, \ee^{- \overline{V}(y)}$.

We now assume furthermore $\psi'(1)>0$, so that $\theta<1$. We choose $\delta\in (0, \, {1\over \theta} -1)$. Since $(1+\delta)\theta<1$, we have
$$
\varrho_n^{ (1+\delta)\theta}
\le
\sum_{j=1}^n \sum_{|y|=j} {\bf 1}_{E_\delta(y)} \,
\ee^{-(1+\delta) \theta \overline{V}(y)}.
$$

\noindent Consider the branching random walk $\widetilde{V}(x):= \theta V(x)$ for any $x$. If we define $\widetilde \psi(t):= \log \E[ \sum_{|x|=1} \ee^{- t\widetilde V(x)}]$, then $\widetilde\psi(1)=\widetilde \psi'(1)=0$. We apply formula (\ref{change-proba}) to $(\widetilde V (x))$, and obtain a centered one-dimensional random walk $(\widetilde S_i,\, 0\le i\le n)$ with $\widetilde \sigma^2 := \E(\widetilde S_1^2)= \E[ \sum_{|x|=1} \theta^2 V(x)^2 \ee^{-\theta V(x)}]$ such that for $1\le j\le n$ (writing $\widehat{S}_i := \max_{1\le k\le i} \widetilde{S}_k$, $\forall i$),
\begin{eqnarray*}
    \E\Big(\sum_{|y|=j} {\bf 1}_{E_\delta(y)} \,
    \ee^{-(1+\delta) \theta \overline{V}(y)} \Big)
 &=& \E\Big(
     \ee^{ \widetilde S_j - (1+\delta) \widehat{S}_j} \,
     {\bf 1}_{ \{ ( 1+\delta) \widehat{S}_i
     - \widetilde{S}_i < a  (n-i)^{1/3}, \,
     \forall i<j, \,
     (1+\delta)\widehat{S}_j -\widetilde{S}_j
     \ge a(n-j)^{1/3} \} }\Big)
     \\
 &\le & \ee^{-a(n-j)^{1/3}}\,
    \P\Big( ( 1+\delta) \widehat{S}_i
    - \widetilde{S}_i < a (n-i)^{1/3}, \, \forall i< j \Big).
\end{eqnarray*}

\noindent It follows that
$$
\E (\varrho_n^{ (1+\delta)\theta})
\le
\sum_{j=1}^n \ee^{-a(n-j)^{1/3}}\,
\P\Big( ( 1+\delta) \widehat{S}_i
- \widetilde{S}_i < a (n-i)^{1/3}, \, \forall i< j \Big).
$$

\noindent We choose $a:= ({3\pi^2 \widetilde \sigma^2 \over 2
})^{1/3} =  \theta \alpha_\theta^{1/3}  $. Applying Corollary
\ref{c:tech} (ii) to $(\widetilde{S}_i)$, we get $\E (\varrho_n^{
(1+\delta)\theta}) \le \ee^{ - (a+o(1)) n^{1/3}}$, for $n\to
\infty$. By Chebyshev's inequality and the Borel--Cantelli lemma,
$\P$-almost surely for $n\to \infty$, $\varrho_n^{ (1+\delta)\theta}
\le \ee^{ - (a+o(1)) n^{1/3}}$. Since $\delta$ can be arbitrarily
small, this implies (\ref{upperrhon}), and completes the proof of
Theorem \ref{t:main2}.\hfill$\Box$

\section{Proof of Theorem \ref{t:main}: upper bound}
\label{s:ub}

$\phantom{aob}$We prove that if $\psi(1)=\psi'(1)=0$, then\footnote{On the set of extinction, the upper bound is, in fact, trivially true.}
\begin{equation}
    \limsup_{n\to \infty} \,
    {\max_{0\le k\le n}|X_k| \over (\log n)^3}
    \le {8\over 3\pi^2 \sigma^2} \, ,
    \qquad\hbox{\rm $\p$-a.s.},
    \label{ub}
\end{equation}

\noindent where $\sigma^2 :=  \E \{ \sum_{|x|=1} V(x)^2 \ee^{- V(x)} \}$.

Let, for any $n\ge 1$,
\begin{equation}
    \beta_n
    :=
    P_\omega\{ \tau_n <
    T_{{\buildrel \leftarrow \over \varnothing}} \},
    \label{beta}
\end{equation}

\noindent where $\tau_n:= \inf\{ i\ge 1: \, |X_i|=n\}$ is as before the first time that the walk reaches the $n$-th generation, whereas $T_{{\buildrel \leftarrow \over \varnothing}} := \inf\{ i\ge 0: X_i = {\buildrel \leftarrow \over \varnothing}\}$ is the first time that the walk hits ${\buildrel \leftarrow \over \varnothing}$. There is a simple relation between $\beta_n$ and $\varrho_n := P_\omega \{ \tau_n < \tau_0 \}$, as stated in the following lemma. We mention that no condition on $\psi$ is in force for the lemma.

\begin{lemma}
 \label{l:beta-rho}
 Assume that the walk $(X_n)$ is recurrent. We have, for all $n\ge 1$,
 \begin{equation}
     \varrho_n
     \le
     \beta_n
     \le {\varrho_n \over
     \omega(\varnothing, \, {\buildrel \leftarrow \over \varnothing})}.
     \label{beta-rho}
 \end{equation}

\end{lemma}

\noindent {\it Proof of Lemma \ref{l:beta-rho}.} The first inequality in (\ref{beta-rho}) is trivial. Let us prove the second. Let $T^{(0)}_\varnothing:=0$ and $T^{(k)}_\varnothing := \inf \{ i> T^{(k-1)}_\varnothing : \; X_i=\varnothing\}$ (for $k\ge 1$). In words, $T^{(k)}_\varnothing$ is the $k$-th return time to the root $\varnothing$. [Thus $T^{(1)}_\varnothing = \tau_0$.] Since the walk is recurrent, each $T^{(k)}_\varnothing$ is well-defined.

Recall that $\beta_n$ represents the probability that, starting from the root, the walk visits generation $n$ before hitting ${\buildrel \leftarrow \over \varnothing}$. By considering the number of returns to $\varnothing$ (which can be 0) by the walk before visiting generation $n$, we have
$$
\beta_n
=
P_\omega\{ \tau_n < T_{{\buildrel \leftarrow \over \varnothing}} \}
=
\sum_{k=0}^\infty
P_\omega \Big\{ T^{(0)}_\varnothing < T^{(1)}_\varnothing < \cdots <
T^{(k)}_\varnothing < \tau_n < T^{(k+1)}_\varnothing, \;
\tau_n < T_{{\buildrel \leftarrow \over \varnothing}} \Big\} .
$$

\noindent Applying the strong Markov property successively at $T^{(k)}_\varnothing$, $\cdots$, $T^{(1)}_\varnothing$, we see that the probability on the right-hand side equals $[P_\omega\{ T^{(1)}_\varnothing < (\tau_n \wedge T_{{\buildrel \leftarrow \over \varnothing}}) \}]^k \, P_\omega \{ \tau_n< T^{(1)}_\varnothing\}$ (notation: $u\wedge v := \min\{ u, \, v\}$). Therefore
$$
\beta_n
=
{P_\omega \{ \tau_n< T^{(1)}_\varnothing\} \over 1-P_\omega\{ T^{(1)}_\varnothing < (\tau_n \wedge T_{{\buildrel \leftarrow \over \varnothing}}) \} }
=
{\varrho_n\over 1-P_\omega\{ \tau_0< (\tau_n \wedge T_{{\buildrel \leftarrow \over \varnothing}}) \} }.
$$

\noindent Since $1-P_\omega\{ \tau_0 < (\tau_n \wedge T_{{\buildrel \leftarrow \over \varnothing}}) \} \ge 1-P_\omega\{ \tau_0 <T_{{\buildrel \leftarrow \over \varnothing}} \} = \omega(\varnothing, \, {\buildrel \leftarrow \over \varnothing})$, this yields the lemma.\hfill$\Box$

\bigskip

We now turn to the proof of (\ref{ub}). Assume $\psi(1)=\psi'(1)=0$. We claim that it suffices to prove that
\begin{equation}
    \limsup_{n\to \infty} \,
    {1\over n^{1/3}} \log \E(\beta_n)
    \le
    - \Big({3\pi^2 \sigma^2 \over 8}\Big)^{1/3} .
    \label{upp-rhon}
\end{equation}

\noindent Indeed, if (\ref{upp-rhon}) holds, then by Chebyshev's inequality and the Borel--Cantelli lemma, for any $\varepsilon>0$ and $\P$-almost surely all sufficiently large $n$, $\beta_n \le \exp[ -(1-\varepsilon) ({3\pi^2 \sigma^2 \over 8})^{1/3} n^{1/3}]$, which by Lemma \ref{l:beta-rho} yields $\varrho_n \le \exp[ -(1-\varepsilon) ({3\pi^2 \sigma^2 \over 8})^{1/3} n^{1/3}]$. In view of Fact \ref{f:rhonxk}, we obtain (\ref{ub}).

It remains to prove (\ref{upp-rhon}). Let $a:= ({3\pi^2 \sigma^2 \over 8})^{1/3}$ and $n\ge 1$. By (\ref{rho-ub}) and Lemma \ref{l:beta-rho},
$$
\E(\beta_n) \le \sum_{j=1}^n \E \Big( \sum_{|y|=j} {\bf 1}_{E_0(y)} \, \ee^{-\overline{V}(y)} \Big),
$$

\noindent where
$$
E_0(y) := \Big\{ \overline{V}(y) - V(y) \ge  a (n-j)^{1/3}\Big\} \cap
\bigcap_{i=1}^{j-1} \Big\{ \overline{V}(y_i) - V(y_i) <  a
(n-i)^{1/3} \Big\} .
$$

\noindent Applying (\ref{change-proba}), this leads to (with $\overline{S}_j := \max_{1\le i\le j} S_i$ as before):
\begin{eqnarray*}
    \E(\beta_n)
 &\le& \sum_{j=1}^n \E \Big\{ \ee^{S_j}
    {\bf 1}_{\{ \overline{S}_j -S_j \ge a (n-j)^{1/3}, \;
    \overline{S}_i -S_i < a (n-i)^{1/3}, \; \forall i<j \}}
    \ee^{- \overline{S}_j} \Big\}
    \\
 &\le&  \sum_{j=1}^n  \ee^{-a (n-j)^{1/3}} \,
    \P\Big\{  \overline{S}_i -S_i < a (n-i)^{1/3}, \;
    \forall i<j \Big\} ,
\end{eqnarray*}


\noindent which, according to Corollary \ref{c:tech} (i), is bounded by $\exp[ - (1+o(1)) ({3\pi^2 \sigma^2 \over 8})^{1/3}  n^{1/3}]$ for $n\to \infty$. This yields (\ref{upp-rhon}).\hfill$\Box$

\section{Proof of Theorem \ref{t:main}: lower bound}
\label{s:lb}

$\phantom{aob}$We start by recalling a spinal decomposition for the branching random walk $(V(x))$.  This decomposition has been used in the literature by many authors in various forms, going back at least to Kahane and Peyri\`ere~\cite{kahane-peyriere}. The material in this paragraph is borrowed from Lyons, Pemantle and Peres~\cite{lyons-pemantle-peres} and Lyons~\cite{lyons}. The starting point is to a change-of-probabilities technique on the space of trees; we refer to the aforementioned  references for more precision.

Assume $\psi(1)=0$, i.e.,  $\E\{ \sum_{|x|=1} \ee^{- V(x)}\}=1$. Let
$$
W_n:= \sum_{|x|=n} \ee^{- V(x)}, \qquad n\ge 0.
$$

\noindent Clearly, $(W_n)$ is a martingale with respect to the filtration $(\mathscr{F}_n)$, where $\mathscr{F}_n$ is the sigma-algebra generated by the branching random walk in the first $n$ generations.

By Kolmogorov's extension theorem, there exists a probability $\Q$
on $\mathscr{F}_\infty$ (the sigma-algebra generated by the
branching random walk) such that for any $n$,
\begin{equation}
    \Q_{| _{ \mathscr{F}_n}} =
    W_n \bullet \P_{| _{ \mathscr{F}_n}} ,
    \label{Q}
\end{equation}

\noindent i.e., $\Q(A) = \E(W_n \, {\bf 1}_A)$, $\forall A\in \mathscr{F}_n$. The law of the branching random walk under the new
probability $\Q$ is called the law of a {\it size-biased branching
random walk}. It is clear that the size-biased branching random walk
survives with probability one.

There is a one-to-one correspondence between a branching random
walk and a  marked tree.  On the  enlarged  probability space formed
by  marked trees with distinguished rays,      we may construct a
probability $\Q$ satisfying (\ref{Q}), and an infinite ray $\{ w_0=
\varnothing, w_1, ..., w_n, ..\}$ such that for any $n\ge 1$, ${\buildrel \leftarrow \over w_n}= w_{n-1}$ (recalling that ${\buildrel \leftarrow \over x}$ is the parent of $x$) and
\begin{equation}
    \Q\Big\{ w_n = x \, \Big| \,
    \mathscr{F}_n \Big\} =
    {\ee^{-V(x)} \over W_n}, \qquad \forall\,   |x|=n.
    \label{u*}
\end{equation}




For any individual $x\not= \varnothing$, let
$$
\Delta V(x) := V(x) - V({\buildrel \leftarrow \over x}) .
$$

\noindent We write, for $k \ge 1$,
\begin{equation}
    \mathscr{I}_k
    :=
    \bigl\{ x : \, |x| =k, \;
    {\buildrel \leftarrow \over x} = w_{k-1} , \;
    x\not= w_k  \bigr\} .
    \label{Ck}
\end{equation}

\noindent In words, $\mathscr{I}_k $ is the set of children of $w_{k-1} $ except $w_k $, or equivalently, the set of the brothers of $w_k$,
and is possibly empty. Finally, let us introduce the following
sigma-field:
\begin{equation}
    \mathscr{G}_n :=
    \sigma \Big\{
    (\Delta V(x), \; x\in\mathscr{I}_k), \;
    V(w_k), \;  w_k , \; \mathscr{I}_k , \; 1\le k\le n\Big\} .
    \label{Gn}
\end{equation}


The promised spinal decomposition is as follows ($xu$ denoting concatenation of $x$ and $u$). Although it slightly differs from the spinal decomposition presented in Lyons~\cite{lyons}, we feel free to omit the proof.

\begin{proposition}
 \label{p:change-proba}
 Assume $\psi(1)=0$, and fix $n\ge1$.
 Under probability $\mathbf{Q}$,

 {\rm (i)} the random variables
 $(\Delta V(w_k), \; \Delta V(x), \; x\in\mathscr{I}_k)$,
 $1\le k\le n$, are i.i.d.;

 {\rm (ii)} conditionally on $\mathscr{G}_n$, the shifted branching
 random walks $(\{ V(xu) - V(x)\}_{|u|=k}, 0\le k\le n-|x|)$, for
 $x\in\bigcup_{k=1}^n \mathscr{I}_k$, are independent, and have the
 same law as $(\{V(u)\}_{|u|=k}, 0\le k\le n-|x|)$ under $\P$.
\end{proposition}

\begin{figure}[h]
\centering
\scalebox{1} 
{
\begin{pspicture}(0,-2.7)(6.839987,2.72)
\definecolor{color747b}{rgb}{0.6,0.6,0.6}
\psline[linewidth=0.02cm](3.58,2.3)(6.48,1.52)
\psline[linewidth=0.04cm](3.6,2.28)(3.82,0.68)
\psline[linewidth=0.02cm](3.6,2.3)(1.76,1.48)
\psline[linewidth=0.02cm](3.6,2.3)(0.38,1.5)
\pstriangle[linewidth=0.0020,dimen=outer,fillstyle=solid,fillcolor=color747b](0.38,-2.68)(0.76,4.2)
\pstriangle[linewidth=0.0020,dimen=outer,fillstyle=solid,fillcolor=color747b](1.77,-2.68)(0.78,4.18)
\pstriangle[linewidth=0.0020,dimen=outer,fillstyle=solid,fillcolor=color747b](6.5,-2.68)(0.68,4.18)
\psline[linewidth=0.02cm](3.8,0.7)(5.16,-0.08)
\psline[linewidth=0.02cm](3.8,0.7)(3.96,-0.1)
\psline[linewidth=0.04cm](3.8,0.68)(3.2,-0.9)
\usefont{T1}{ptm}{m}{n}
\rput(4.0,2.54){\footnotesize $w_0$}
\usefont{T1}{ptm}{m}{n}
\rput(4.34,0.92){\footnotesize $w_1$}
\pstriangle[linewidth=0.0020,dimen=outer,fillstyle=solid,fillcolor=color747b](5.17,-2.68)(0.78,2.6)
\pstriangle[linewidth=0.0020,dimen=outer,fillstyle=solid,fillcolor=color747b](3.97,-2.68)(0.78,2.6)
\usefont{T1}{ptm}{m}{n}
\rput(2.82,-0.8){\footnotesize $w_2$}
\usefont{T1}{ptm}{m}{n}
\rput(2.7,-2.48){\footnotesize $w_n$}
\psline[linewidth=0.04cm,linestyle=dotted,dotsep=0.16cm](3.21,-0.88)(3.19,-2.54)
\usefont{T1}{ptm}{m}{n}
\rput(1.78,-2.34){\footnotesize $\P$}
\usefont{T1}{ptm}{m}{n}
\rput(0.38,-2.32){\footnotesize $\P$}
\usefont{T1}{ptm}{m}{n}
\rput(3.98,-2.28){\footnotesize $\P$}
\usefont{T1}{ptm}{m}{n}
\rput(6.4,-2.28){\footnotesize $\P$}
\usefont{T1}{ptm}{m}{n}
\rput(5.2,-2.3){\footnotesize $\P$}
\end{pspicture}
}

\caption{A \Q-tree}
\end{figure}

We now proceed to (the beginning of) the proof of the lower bound in Theorem \ref{t:main}, of which we recall the statement: under the assumption $\psi(1)=\psi'(1)=0$, we have, on the set of non-extinction,
\begin{equation}
    \liminf_{n\to \infty} \,
    {\max_{0\le k\le n}|X_k| \over (\log n)^3}
    \ge
    {4\over \alpha} = {8\over 3\pi^2 \sigma^2} \, ,
    \qquad\hbox{\rm $\p$-a.s.},
    \label{lb}
\end{equation}

\noindent where $\sigma^2 := \E \{ \sum_{|x|=1} V(x)^2 \ee^{- V(x)} \}$.

Let $\beta_n := P_\omega\{ \tau_n < T_{{\buildrel \leftarrow \over \varnothing}} \}$ be as in (\ref{beta}), where $\tau_n= \inf\{ i\ge 1: \, |X_i|=n\}$, and $T_{{\buildrel \leftarrow \over \varnothing}} = \inf\{ i\ge 0: X_i = {\buildrel \leftarrow \over \varnothing}\}$. We claim that it suffices to prove that
\begin{equation}
    \liminf_{n \to\infty} \, {1\over n^{1/3}} \log \E(\beta_n)
    \ge
    - \Big( {3\pi^2\sigma^2\over 8}\Big)^{1/3} .
    \label{low-epn}
\end{equation}

It is indeed easy to check that (\ref{low-epn}) implies (\ref{lb}): Let $\mathscr{S} := \{ \hbox{the system survives}\}$, $\mathscr{S}_n := \{ \hbox{the system survives at least until generation $n$}\}$. Clearly $\mathscr{S} \subset \mathscr{S}_n$ for any $n$. Recall that there exists (see \cite{yzpolymer}, p.~755) a constant $c>0$ such that for all large $n$,
$$
\P(W_n < n^{-c} \, | \, \mathscr{S}_n) \le n^{-2}.
$$

\noindent On the other hand, we have (see \cite{yztree}, p.~543,
Remark; the result therein states for the regular tree, but the
same proof by convexity obviously holds in the general case)
$$
\E \Big( \ee^{- t {\beta_n \over \E(\beta_n)}}\Big)
\le
\E (\ee^{-t W_n}), \qquad t\ge 0.
$$

\noindent Since $\beta_n=0=W_n$ on $\mathscr{S}_n^c$, it is equivalent to say that $\E ( \ee^{- t {\beta_n \over \E(\beta_n)}} \, | \, \mathscr{S}_n) \le \E (\ee^{-t W_n} \, | \, \mathscr{S}_n)$. Therefore, for any $\varepsilon>0$ and all sufficiently large $n$,
$$
\P \Big( {\beta_n \over \E(\beta_n)} < \ee^{-\varepsilon n^{1/3}} \, \Big| \, \mathscr{S}_n \Big)
\le
\ee^1\,\E \Big( \ee^{- \ee^{\varepsilon n^{1/3}} W_n}\, \Big| \, \mathscr{S}_n \Big)
\le
n^{-2} \ee + \ee^{- n^{-c} \ee^{\varepsilon n^{1/3}} } .
$$

\noindent Since $\mathscr{S} \subset \mathscr{S}_n$, this implies $\sum_n \P ( {\beta_n \over \E(\beta_n)} < \ee^{-\varepsilon n^{1/3}} \, | \, \mathscr{S} ) \le {1\over \P(\mathscr{S})}\sum_n \P ( {\beta_n \over \E(\beta_n)} < \ee^{-\varepsilon n^{1/3}} \, | \, \mathscr{S}_n) <\infty$. If (\ref{low-epn}) holds, then by the Borel--Cantelli lemma, on the set $\mathscr{S}$, $\P$-almost surely for all sufficiently large $n$, $\beta_n \ge \ee^{- \varepsilon n^{1/3}} \E(\beta_n) \ge \exp\{ -[2\varepsilon+ ( {3\pi^2\sigma^2\over 8})^{1/3}] n^{1/3}\}$, and thus $\varrho_n \ge \omega(\varnothing, \, {\buildrel \leftarrow \over \varnothing}) \exp\{ -[2\varepsilon+ ( {3\pi^2\sigma^2\over 8})^{1/3}] n^{1/3}\}$ (Lemma \ref{l:beta-rho}). In view of Fact \ref{f:rhonxk}, we obtain (\ref{lb}), the lower bound in Theorem \ref{t:main}.

The rest of the section is devoted to the proof of (\ref{low-epn}). Let as before $\varrho_n := P_\omega \{ \tau_n < \tau_0\}$. Since $\beta_n \ge \varrho_n$ (Lemma \ref{l:beta-rho}), we only need to bound $\E(\varrho_n)$ from below.

For any vertex $x$, let $P_\omega^x$ be the (quenched) probability such that $P_\omega^x \{ X_0=x\} =1$. We first prove a formula for $\varrho_n$ without the assumption $\psi(1)=\psi'(1)=0$. We mention that if $|x|=n$, then under $P_\omega^x$, $\tau_n$ is the first {\it return} time to generation $n$.

\begin{lemma}
\label{l:rho}
 Assume that the walk $(X_n)$ is recurrent.
 For any $n\ge 1$, we have
 $$
 \varrho_n =
 \omega(\varnothing, \,
 {\buildrel \leftarrow \over \varnothing})
 \sum_{|x|=n}
 {\ee^{-V(x)} \over
 \omega(x, \, {\buildrel \leftarrow \over x})}
 P_\omega^x \{ \tau_n> \tau_0\} .
 $$

\end{lemma}

\noindent {\it Proof of Lemma \ref{l:rho}.} The beginning of the proof uses a similar idea as in the proof of Lemma \ref{l:beta-rho}, except that instead of considering the number of returns to $\varnothing$ before hitting generation $n$, we consider the last site at generation $n$ visited by the walk during an excursion. More precisely, for any $x$ with $|x| \ge 1$, let $T^{0)}_x :=0$ and
$T^{(k)}_x := \inf \{ i> T^{(k-1)}_x: \; X_i=x\}$ (for $k\ge 1$). In words, $T^{(k)}_x$ is the time of the $k$-th visit at $x$.

Recall that $\varrho_n$ is the (quenched) probability that during an excursion away from the root $\varnothing$, the walk hits generation $n$. By considering the last site at generation $n$ visited by the walk during the excursion, we have
\begin{eqnarray*}
    \varrho_n
 &=& \sum_{|x|=n} \sum_{k=1}^\infty
    P_\omega \Big\{ T^{(k)}_x < \tau_0
    < T^{(k+1)}_x, \;
    \max_{T^{(k)}_x <i \le \tau_0} |X_i|<n
    \Big\}
    \\
 &=& \sum_{|x|=n} \sum_{k=1}^\infty
    P_\omega \Big\{ T^{(k)}_x < \tau_0, \;
    \max_{T^{(k)}_x <i \le \tau_0} |X_i|<n
    \Big\} .
\end{eqnarray*}

\noindent Applying the strong Markov property at $T^{(k)}_x$, we
see that the probability on the right-hand side equals $P_\omega\{
T^{(k)}_x < \tau_0\} \, P_\omega^x \{ \tau_n> \tau_0\}$. Therefore,
$$
\varrho_n = \sum_{|x|=n} P_\omega^x \{ \tau_n> \tau_0\}
\sum_{k=1}^\infty P_\omega\{ T^{(k)}_x < \tau_0\} = \sum_{|x|=n}
P_\omega^x \{ \tau_n> \tau_0\} E_\omega\Big( \sum_{i=0}^{\tau_0-1}
{\bf 1}_{ \{ X_i=x\} } \Big).
$$

\noindent  $E_\omega ( \sum_{i=0}^{\tau_0-1} {\bf 1}_{ \{
X_i=x\} })$, is the expected number of visits at site $x$ in an
excursion, and can therefore be explicitly computed.  Indeed one can easily check that, as a function of $x$, it is invariant with respect to the transition matrix $\omega(x,y)$. In the particular setting of Markov chains on trees any invariant measure can be computed, using an easy recurrence. One gets that all the invariant measures are proportional to $\pi(x) := {1\over \omega(x, \, {\buildrel
\leftarrow \over x})} \, \ee^{-V(x)}$, for $x\not= {\buildrel
\leftarrow \over \varnothing}$ (Note that this formula is also valid for $x=\varnothing,$ because of the consistent definition of $\omega(\varnothing, \, {\buildrel
\leftarrow \over \varnothing})$). Therefore, there exists $0<c(\omega)<\infty$ such that
$$
E_\omega\Big( \sum_{i=0}^{\tau_0-1} {\bf 1}_{ \{ X_i=x\} } \Big) =
{c(\omega) \over \omega(x, \, {\buildrel \leftarrow \over x})} \,
\ee^{-V(x)} .
$$

\noindent To determine the value of $c(\omega)$, we take $x:=
\varnothing$, to see that $c(\omega) = \omega(\varnothing, \, {\buildrel
\leftarrow \over \varnothing})$. This yields the lemma.\hfill$\Box$

\bigskip

Assume $\psi(1)=0$. We make use of the size-biased branching random walk, and work under the new probability $\Q$. Recall the definitions of $\Q$ and $w_n$ from (\ref{Q}) and (\ref{u*}), respectively. By Lemma \ref{l:rho},
$$
\E(\varrho_n) = \E_\Q \Big\{ {\omega(\varnothing, \, {\buildrel
\leftarrow \over \varnothing}) \over \omega(w_n, \, w_{n-1}) }
P_\omega^{w_n} \{ \tau_n> \tau_0\} \Big\} .
$$

\noindent We observe that
$$
P_\omega^{w_n} \{ \tau_n> \tau_0\} = \prod_{j=1}^n P_\omega^{w_j} \{
\tau_n> T(w_{j-1})\} =: \prod_{j=1}^n Y_j.
$$

\noindent Obviously, $Y_n = \omega(w_n, \, w_{n-1})$, $Y_{n-1} = \omega(w_{n-1}, \, w_{n-2})$.

Let $j\le n-2$. By the Markov property, $Y_j = \omega(w_j, \, w_{j-1}) + \sum_{x: \; {\buildrel \leftarrow \over x} = w_j} \omega(w_j, \, x)
P_\omega^x \{ \tau_n> T(w_{j-1})\}$, whereas by the strong Markov
property, $P_\omega^x \{ \tau_n> T(w_{j-1})\} = P_\omega^x \{
\tau_n> T(w_j)\}\, Y_j$ for all $x$ such that ${\buildrel \leftarrow
\over x} = w_j$. Accordingly,
$$
Y_j
=
{\omega(w_j, \, w_{j-1})\over
    1- \sum_{x: \; {\buildrel \leftarrow \over x} = w_j}
    \omega(w_j, \, x) \, P_\omega^x \{ \tau_n >T(w_j)\} }
=
{1\over
    1 + \sum_{x: \; {\buildrel \leftarrow \over x} = w_j}
    B(x) P_\omega^x \{ \tau_n< T(w_j)\} } ,
$$

\noindent where
$$
B(x) := \ee^{-[V(x) - V({\buildrel \leftarrow \over x})]} =
{\omega({\buildrel \leftarrow \over x}, \, x) \over
\omega({\buildrel \leftarrow \over x}, \, {\buildrel \Leftarrow
\over x})} .
$$

\noindent So, if we write
$$
\xi_j := \sum_{x: \;
    {\buildrel \leftarrow \over x} = w_j, \;
    x\not= w_{j+1} }
    B(x) P_\omega^x \{ \tau_n< T(w_j)\} , \qquad 1\le j\le n-2,
$$

\noindent then $Y_j = {1\over 1+\xi_j + (1-Y_{j+1}) B(w_{j+1})}$, $1\le j\le n-2$, and $\E(\varrho_n) = \E_\Q \{ {\omega(\varnothing, \, {\buildrel \leftarrow \over \varnothing}) \over \omega(w_n, \, w_{n-1}) } \prod_{j=1}^n Y_j \}= \E_\Q \{ \omega(\varnothing, \, {\buildrel \leftarrow \over \varnothing}) \prod_{j=1}^{n-1} Y_j \}$.

Let $\mathscr{G}_n$ be the sigma-algebra generated by the first $n$ generations of the spine (see (\ref{Gn})). By Proposition \ref{p:change-proba}, under $\Q$, the random variables $\xi_1$, $\cdots$, $\xi_{n-1}$ are conditionally independent given $\mathscr{G}_n$. Moreover, for any $1\le j \le n-2$,
\begin{equation}
    \E_\Q(\xi_j  \, |\, \mathscr{G}_n)
    =
    \sum_{x: \;
    {\buildrel \leftarrow \over x} = w_j, \; x\not= w_{j+1} }
    B(x) \, \E(\beta_{n-1-j})
    \le
    {\E(\beta_{n-1-j} ) \over \omega(w_j, w_{j-1})}  .
    \label{E(xi)<}
\end{equation}

We now provide a lower bound for $\E(\varrho_n)$,  by replacing
$(Y_j)_{1\le j\le n-1}$ by a new collection of random variables,
denoted by $(Z_j)_{1\le j\le n-1}$ and defined as follows: $Z_{n-1} := Y_{n-1} = \omega(w_{n-1}, \, w_{n-2})$ and for $1\le j\le n-2$,
\begin{equation}
    Z_j
    :=
    {1\over 1+ \E_\Q(\xi_j \, | \, \mathscr{G}_n)
    + (1-Z_{j+1})B(w_{j+1})}.
    \label{zjfj}
\end{equation}

\noindent Since $Z_{n-1}$, $B(w_{n-1})$, $B(w_{n-2})$, $\cdots$, $B(w_2)$ are $\mathscr{G}_n$-measurable, it follows by backwards induction on $j$ that each $Z_j$, for $1\le j\le n-1$, is $\mathscr{G}_n$-measurable.

\begin{lemma}
\label{l:E(rho|G)>}
 Assume $\psi(1)=0$.
 For any $n\ge 3$, we have
 $$
 \E_\Q \Big\{
 \prod_{j=1}^{n-1} Y_j \, \Big| \, \mathscr{G}_n
 \Big\}
 \ge \prod_{j=1}^{n-1} Z_j ,
 \qquad\hbox{\rm $\Q$-a.s.}
 $$

\end{lemma}

\noindent {\it Proof of Lemma \ref{l:E(rho|G)>}.} For any $c\in [0,
\, 1]$ and $a:= (a_1, \cdots, a_{n-1}) \in \r_+^{n-1}$, we define
$F_{n-1}^{c,a} (u_{n-1}) := c$, $u_{n-1} \in \r_+$, and for $1\le
j\le n-2$,
$$
F_j^{c,a} (u_j, \cdots, u_{n-2}) := {1\over 1+ u_j +
a_{j+1} [1- F_{j+1}^{c,a} (u_{j+1}, \cdots, u_{n-2})]}, \quad
(u_j, \cdots, u_{n-2}) \in \r_+^{n-j-1}.
$$

\noindent Then by backwards induction on $j$, we have, for $1\le j\le n-1$,
$$
Y_j =
F_j^{Y_{n-1}, B(w)}(\xi_j, \cdots, \xi_{n-2}),
\qquad
Z_j =
F_j^{Z_{n-1}, B(w)}(\E_\Q(\xi_j \, | \, \mathscr{G}_n), \cdots, \E_\Q(\xi_{n-2} \, | \, \mathscr{G}_n)),
$$

\noindent where $B(w) := (B(w_1), \cdots, B(w_{n-1}))$. Note that both $Y_{n-1}$ and $B(w)$ are $\mathscr{G}_n$-measurable.

Recall that (under $\Q$) $\xi_1$, $\cdots$, $\xi_{n-2}$ are conditionally independent given $\mathscr{G}_n$. By Jensen's inequality, if $\Phi: \, \r^{n-2}_+ \to \r$ is coordinate-wise convex, then $\E_\Q \{ \Phi(\xi_1, \, \cdots, \xi_{n-2})\, | \, \mathscr{G}_n)\} \ge \Phi(\E_\Q(\xi_1 \, | \, \mathscr{G}_n), \cdots, \E_\Q(\xi_{n-2}\, | \, \mathscr{G}_n))$, $\Q$-a.s. So we only need to show that for any $c\in [0, \, 1]$ and $a \in \r_+^{n-1}$, $(u_1, \, \cdots, u_{n-2})\mapsto \prod_{j=1}^{n-1} F_j^{c,a}(u_j, \cdots, u_{n-2})$ as a function on $\r_+^{n-2}$, is convex in each of $u_i$.

Since the product of non-negative, coordinate-wise non-increasing,
coordinate-wise convex functions is still (non-negative,
coordinate-wise non-increasing, and) coordinate-wise convex, we only
have to check that for any $j\le n-2$, the function $(u_j, \, \cdots,
u_{n-2}) \mapsto F_j^{c,a}(u_j, \cdots, u_{n-2})$ is
non-negative (which is obvious), coordinate-wise non-increasing, and
coordinate-wise convex. We prove it by induction on $j$.

By definition, $F_{n-2}^{c,a}(u_{n-2}) = [1+ u_{n-2} +
(1-c) a_{n-1}]^{-1}$, which is obviously non-increasing and convex in
$u_{n-2}$.

Assume that for $1\le j\le n-3$, $(u_{j+1}, \, \cdots, u_{n-2})
\mapsto F_{j+1}^{c,a}(u_{j+1}, \cdots, u_{n-2})$ is
coordinate-wise non-increasing and coordinate-wise convex. Since
$$
F_j^{c,a} (u_j, \cdots, u_{n-2}) = {1\over 1+ u_j +
a_{j+1} [1- F_{j+1}^{c,a} (u_{j+1}, \cdots, u_{n-2})]} ,
$$

\noindent $F_j^{c,a}$ is non-increasing and convex in each of $u_i$ (for $j\le i\le n-2$): the monotonicity is obvious, whereas the convexity follows from the fact that $y\mapsto {1\over 1+ u_j + (1-y)a_{j+1}}$ is convex and non-decreasing on $[0, \, 1]$ and that $f\circ g$ is convex if $f$ is convex and non-decreasing while $g$
is convex.\hfill$\Box$

\bigskip

Recall that $\E(\varrho_n) = \E_\Q \{ \omega(\varnothing, \,
{\buildrel \leftarrow \over \varnothing}) \prod_{j=1}^{n-1} Y_j \}$. Since $\omega(\varnothing, \, {\buildrel \leftarrow \over \varnothing})$ is $\mathscr{G}_n$-measurable, it follows from Lemma
\ref{l:E(rho|G)>} that
\begin{equation}
    \E(\varrho_n)
    \ge
    \E_\Q \Big\{
    \omega(\varnothing, \, {\buildrel \leftarrow \over \varnothing})
    \prod_{j=1}^{n-1} Z_j \Big\}.
    \label{erho1}
\end{equation}

We now give a lower bound for $\prod_{j=1}^{n-1} Z_j$ by means of a deterministic lemma. The proof of the lemma is in the Appendix.

\begin{lemma}
\label{markov-zj}
 Let $n>  k\ge 2$. Let $b_{j+1}>0$ and $r_j \ge 0$ for all
 $0\le j < n$. Define $(z_j)_{1\le j\le n}$ by $z_n=0$ and
 $$
 z_j
 :=
 {1\over 1+  r_j + b_{j+1} (1-z_{j+1}) }, \qquad 1\le j\le n-1.
 $$
 Let $v(0):=0$ and $v(j):= -\sum_{i=1}^j \log b_i$, $1\le j \le n$.
 For any $m_0=0< m_1< ...< m_k=n-1$, we have
 $$
 \prod_{j=1}^{n-1} z_j
 \ge
 {2^{-k} \over \prod_{i=1}^k (m_i- m_{i-1})}
 \, \exp \Big\{  -
 \sum_{i=1}^k \Big( \lambda_i +
 (m_i-m_{i-1})^2\, r^{(i)} \, \ee^{v^*_i}\Big)\Big\},
 $$
 where for $1\le i\le k$ (with $y^+ := \max\{y, \, 0\}$ for $y\in \r$),
 \begin{eqnarray*}
     r^{(i)}
  &:=& \max_{m_{i-1} < j \le m_i} r_j ,
     \\
     \lambda_i
  &:=& \max_{m_{i-1} <j \le m_i} (v(j) - v(m_i))
     + (v(m_i) - v(1+m_i))^+ ,
     \\
     v^*_i
  &:=& \max_{m_{i-1}< j\le \ell \le m_i} (v(j)-v(\ell)).
 \end{eqnarray*}
\end{lemma}

We continue with the proof of the lower bound in Theorem \ref{t:main}. Recall from (\ref{erho1}) that $\E(\varrho_n) \ge \E_\Q \{ \omega(\varnothing, \, {\buildrel \leftarrow \over \varnothing}) \prod_{j=1}^{n-1} Z_j \}$.

Let $k\ge 2$ and $m_0:=0<m_1<m_2<...<m_k=n-1$. Taking $b_{j+1}= B(w_{j+1})$ and $r_j := \E_\Q(\xi_j  \, |\, \mathscr{G}_n)$, we note from (\ref{zjfj}) that we may take the choice of $z_j=Z_j$ in Lemma \ref{markov-zj}. Applying this lemma, and arguing that $\prod_{i=1}^k (m_i - m_{i-1}) \le \prod_{i=1}^k n = n^k$, we find that
\begin{equation}
    \E(\varrho_n)
    \ge
    {1\over (2n)^k} \,
    \E_\Q\Big(
    \omega(\varnothing, \, {\buildrel \leftarrow \over \varnothing}) \,
    \ee^{- \sum_{i=1}^k \Lambda_i
    - \sum_{i=1}^k (m_i - m_{i-1})^2  \, r^{(i)} \,
    \ee^{S^*_i}}\Big),
    \label{fineho}
\end{equation}

\noindent where, for any $1\le i \le k$,
\begin{eqnarray*}
    r^{(i)}
 &:=& \max_{m_{i-1} < j \le m_i} \E_\Q(\xi_j  \, |\, \mathscr{G}_n)
    \le \max_{m_{i-1} < j \le m_i}
    {\E(\beta_{n-1-j} ) \over \omega(w_j, w_{j-1})},
    \\
    \Lambda_i
 &:=& \max_{m_{i-1} < j \le m_i} (S_j - S_{m_i})
    + ( S_{m_i}- S_{1+m_i})^+,
    \\
    S^*_i
 &:=& \max_{m_{i-1}<j\le \ell \le m_i} (S_j -S_\ell),
\end{eqnarray*}

\noindent with $S_j := V(w_j)$, $0\le j\le n$. [In the inequality for $r^{(i)}$, we used (\ref{E(xi)<}).]


We choose: $\chi := {1\over 100}$, $k:= \lfloor n^{1-\chi\over 3}\rfloor$, $m_0:=0$, $m_i:= n- (k-i)^3 \lfloor n^\chi\rfloor$ for $1\le i \le k-1$, and $m_k:=n-1$.

Let $c>1$ be a constant sufficiently large such that $\Q\{ S_2\ge S_1, \; \omega(w_1, \, \varnothing) \ge {1\over c}\} > {1\over c}$. Let
$$
E_n^{(1)}
:=
\bigcap_{j= m_{k-1} +1}^{m_k}
\left\{ S_{j+1} \ge S_j, \;
    \omega(w_j, w_{j-1}) \ge {1\over c}\right\} .
$$

\noindent On $E_n^{(1)}$, we have $\Lambda_k \le 0$, $r^{(k)} \le c$, and $S^*_k=0$, whereas by definition, $m_k-m_{k-1} = \lfloor n^\chi\rfloor -1 \le n^\chi$. Therefore, by (\ref{fineho}),
\begin{eqnarray*}
    \E(\varrho_n)
 &\ge&
    {\ee^{- c\, n^{2\chi} } \over (2n)^k}
    \, \E_\Q\Big(
    \omega(\varnothing, \, {\buildrel \leftarrow \over \varnothing}) \,
    \ee^{- \sum_{i=1}^{k-1} \Lambda_i -
    \sum_{i=1}^{k-1} (m_i - m_{i-1})^2  \, r^{(i)} \,
    \ee^{S^*_i}} \,
    {\bf 1}_{E_n^{(1)}} \Big)
    \\
 &=& {\ee^{- c\, n^{2\chi} } \over (2n)^k}
    \, \E_\Q\Big(
    \omega(\varnothing, \, {\buildrel \leftarrow \over \varnothing}) \,
    \ee^{- \sum_{i=1}^{k-1} \Lambda_i -
    \sum_{i=1}^{k-1} (m_i - m_{i-1})^2  \, r^{(i)} \,
    \ee^{S^*_i}} \Big) \,
    \Q(E_n^{(1)}),
\end{eqnarray*}

\noindent the last identity being a consequence of the fact (notation: $w_{-1} := {\buildrel \leftarrow \over \varnothing}$) that under $\Q$, $(S_j -S_{j-1}, \, \omega(w_{j-1}, w_{j-2}))$, for $j\ge 1$, are independent (they are i.i.d.\ for $j\ge 2$). By the definition of $c$, $\Q( E_n^{(1)} ) = [ \Q\{ S_2\ge S_1, \; \omega(w_1, \, \varnothing) \ge {1\over c}\} ]^{m_k-m_{k-1}} \ge ({1\over c})^{m_k-m_{k-1}} = ({1\over c})^{\lfloor n^\chi\rfloor -1}$. Hence,
$$
\E(\varrho_n)
    \ge
    {\ee^{- c\, n^{2\chi} }
    \over (2n)^k\, c^{\lfloor n^\chi\rfloor -1}}
    \, \E_\Q\Big(
    \omega(\varnothing, \, {\buildrel \leftarrow \over \varnothing}) \,
    \ee^{- \sum_{i=1}^{k-1} \Lambda_i -
    \sum_{i=1}^{k-1} (m_i - m_{i-1})^2  \, r^{(i)} \,
    \ee^{S^*_i}} \Big) .
$$

Let $\varepsilon \in (0, \, {\chi\over 3})$. Write $a_* := ({3\pi^2 \sigma^2\over 8})^{1/3}$. By (\ref{upp-rhon}), there exists some constant $c_1>0$ such that $\E(\beta_i) \le c_1\, \ee^{- (a_*-\varepsilon) (i+1)^{1/3}}$ for all $i\ge1$. Thus
$$
r^{(i)}
\le
{c_1 \, \ee^{- (a_*-\varepsilon) (n- m_i)^{1/3}} \over \min_{m_{i-1} < j \le m_i} \omega(w_j, \, w_{j-1})}, \qquad 1\le i \le k.
$$

\noindent Consider
$$
E_n^{(2)}
:= \left\{ \omega(w_j, w_{j-1}) \ge \ee^{- n^{\varepsilon}}, \,
    \forall\, 1\le j \le m_{k-1} \right\}
    \cap \{\omega(\varnothing, \, {\buildrel \leftarrow \over \varnothing})
    \ge \ee^{- n^{\varepsilon}}\} .
$$

\noindent On $E_n^{(2)}$, we have, for any $1\le i \le k-1$, $r^{(i)} \le c_1 \ee^{ - (a_*- \varepsilon) (n- m_i)^{1/3} + n^\varepsilon}$, whereas $m_i - m_{i-1}\le n$, thus $(m_i - m_{i-1})^2 r^{(i)} \le \ee^{ - (a_*- 2\varepsilon) (n - m_i)^{1/3} }$ (for all sufficiently large $n$; we insist on the fact that $i<k$). Hence
\begin{equation}
    \E(\varrho_n)
    \ge
    {\ee^{- c\, n^{2\chi} -n^\varepsilon}
    \over (2n)^k\, c^{\lfloor n^\chi\rfloor -1}}
    \, \E_\Q\Big( \ee^{- \sum_{i=1}^{k-1} [\Lambda_i +
    \ee^{ S^*_i- (a_*- 2\varepsilon) (n - m_i)^{1/3} }]} \,
    {\bf 1}_{E_n^{(2)}} \Big) .
    \label{fineho2}
\end{equation}

\noindent Let, for $1\le i\le k-1$,
$$
    E_{n,i}^{(3)}
 := \Big\{
    S_i^* < (a_*- 2\varepsilon) (n- m_i)^{1/3} , \,
    \max_{m_{i-1} < j \le m_i} (S_j - S_{m_i})  \le n^\varepsilon, \;
    |S_{1+m_i} - S_{m_i} | \le n^\varepsilon
    \Big\} .
$$

\noindent On the event $E_{n,i}^{(3)}$ (for $1\le i\le k-1$), we have $\Lambda_i \le n^\varepsilon + n^\varepsilon = 2n^\varepsilon$, and, of course, $S_i^* -(a_*- 2\varepsilon) (n - m_i)^{1/3} \le 0$, so that $\Lambda_i + \ee^{ S^*_i- (a_*- 2\varepsilon) (n  - m_i)^{1/3}} \le 2 n^\varepsilon +1 \le 3n^\varepsilon$. Going back to (\ref{fineho2}), we obtain:
\begin{equation}
    \E(\varrho_n)
    \ge
    {\ee^{-c\, n^{2\chi} -n^\varepsilon-3n^\varepsilon(k-1)}
    \over (2n)^k \, c^{\lfloor n^\chi\rfloor -1}}
    \, \Q\Big( E_n^{(2)} \cap
    \bigcap_{i=1}^{k-1} E_{n,i}^{(3)}\Big) .
    \label{fineho3}
\end{equation}

By independence,
\begin{eqnarray*}
    \Q\Big( E_n^{(2)} \cap
    \bigcap_{i=1}^{k-1} E_{n,i}^{(3)}\Big)
 &=&\Q \{\omega(\varnothing,\,{\buildrel \leftarrow \over\varnothing})
    \ge \ee^{- n^{\varepsilon}}\}
    \prod_{i=1}^{k-1} \Q \Big( E_{n,i}^{(3)}, \,
    \min_{m_{i-1} < \ell \le m_i} \omega(w_\ell, \, w_{\ell-1})
    \ge \ee^{- n^\varepsilon} \Big)
    \\
 &\ge& {1\over 2} \, \prod_{i=1}^{k-1} \Q \Big( E_{n,i}^{(3)}, \,
    \min_{m_{i-1} < \ell \le m_i} \omega(w_\ell, \, w_{\ell-1})
    \ge \ee^{- n^\varepsilon} \Big) ,
\end{eqnarray*}

\noindent the last inequality holding for all sufficiently large $n$ (in view of the fact that $\Q \{\omega(\varnothing,\,{\buildrel \leftarrow \over\varnothing}) \ge \ee^{- n^{\varepsilon}}\} \to 1$, $n\to \infty$). By independence again, for any $1\le i\le k-1$,
\begin{eqnarray*}
 &&\Q \Big( E_{n,i}^{(3)}, \,
    \min_{m_{i-1} < \ell \le m_i} \omega(w_\ell, \, w_{\ell-1})
    \ge \ee^{- n^\varepsilon} \Big)
    \\
 &=& \Q \Big( |S_{1+m_i} - S_{m_i} | \le n^\varepsilon, \;
    \omega(w_{m_i}, \, w_{m_i-1})
    \ge \ee^{- n^\varepsilon} \Big) \times
    \\
 && \times \Q \Big( S_i^* < (a_*- 2\varepsilon) (n- m_i)^{1/3} , \,
    \max_{m_{i-1} < j \le m_i} (S_j - S_{m_i})  \le n^\varepsilon, \,
    \\
 && \qquad\qquad
    \min_{m_{i-1} < \ell <m_i} \omega(w_\ell, \, w_{\ell-1})
    \ge \ee^{- n^\varepsilon} \Big)
    \\
 &=& \Q \Big( |S_2 - S_1| \le n^\varepsilon, \;
    \omega(w_1, \, \varnothing)
    \ge \ee^{- n^\varepsilon} \Big) \times
    \Q \Big( F_i(n), \,
    \min_{1\le \ell <\Delta_i} \omega(w_\ell, \, w_{\ell-1})
    \ge \ee^{- n^\varepsilon} \Big) ,
\end{eqnarray*}

\noindent where, for $1\le i\le k-1$,
\begin{eqnarray*}
    \Delta_i
 &:=&m_i-m_{i-1},
    \\
    F_i(n)
 &:=&\Big\{
    \max_{1\le \ell \le \Delta_i}
    ( \overline{S}_\ell - S_\ell) <
    (a_*- 2\varepsilon) (n- m_i)^{1/3}, \;
    \overline{S}_{\Delta_i} -  S_{\Delta_i}  \le n^\varepsilon \Big\},
\end{eqnarray*}

\noindent with $\overline{S}_\ell := \max_{1\le j\le \ell} S_j$ as before. Again, $\Q \{ |S_2 - S_1| \le n^\varepsilon, \; \omega(w_1, \, \varnothing) \ge \ee^{- n^\varepsilon} \}$ is greater than ${1\over 2}$ for large $n$ because it converges to 1. Therefore, for all large $n$,
$$
\Q\Big( E_n^{(2)} \cap
    \bigcap_{i=1}^{k-1} E_{n,i}^{(3)}\Big)
\ge
{1\over 2^k} \, \prod_{i=1}^{k-1} \Q \Big( F_i(n), \,
    \min_{1\le \ell < \Delta_i} \omega(w_\ell, \, w_{\ell-1})
    \ge \ee^{- n^\varepsilon} \Big) .
$$

\noindent To bound the probability expression on the right-hand side, we use the following lemma, which is a uniform version of Proposition \ref{p:sbar}. Its proof is in the Appendix.

\begin{lemma}
\label{L:inter}
 Let $S_i- S_{i-1}$, $i\ge 1$, be i.i.d.\ mean-zero random variables ($S_0:=0$) with $\sigma^2 := \E(S_1^2) \in (0, \, \infty)$.
 For any $\delta >0$, there exist $r_0>1$ and $0<\eta<1$ such that for
 all $r\in [r_0, \, \eta \, \sqrt{n}\,] \cap \z$, for all events
 $A^{(n)}_i$, $1\le i\le n$, satisfying the following two conditions:

 $\bullet$ $(S_i- S_{i-1}, \, A^{(n)}_i)$, for $1\le i \le n$, are
 i.i.d.,

 $\bullet$ $\P( \cap_{i=1}^n A_i^{(n)} ) \ge 1- {\eta \over r}$,

 \noindent we have
 $$
 {\eta \over r} \,
 \ee^{ - (1+\delta) {\pi^2 \sigma^2\over 8} {n\over r^2}}
 \le
 \P\Big( \max_{1\le i \le n}( {\overline{S}_i - S_i}) < r, \,
 \overline{S}_n= S_n, \,
 \bigcap_{i=1}^n A_i^{(n)} \Big)
 \le
 \ee^{ - (1-\delta) {\pi^2\sigma^2\over 8} {n \over r^2}},
 $$
 where $\overline{S}_i:=\max_{1\le j \le i} S_j$.
\end{lemma}

Recall that $(A_1, \cdots, A_N)$ is a random vector distributed as any of $(A_1(x), \cdots, A_{N(x)}(x))$ defined in (\ref{Ai}). Since $\E[ {1\over \omega(\varnothing, \, {\buildrel \leftarrow \over \varnothing})} ]= \E[ 1+ \sum_{i=1}^N A_i] = 1+ \ee^{\psi(1)} =2$, we have, $\Q \{ \omega(w_1, \, \varnothing) < \ee^{-n^\varepsilon} \} = \P \{ \omega(\varnothing, \, {\buildrel \leftarrow \over \varnothing}) < \ee^{-n^\varepsilon} \} \le 2\ee^{-n^\varepsilon}$ (by Markov's inequality). Therefore, for all large $n$ and all $1\le i\le k-1$,
\begin{eqnarray*}
    \Q\Big\{ \min_{1\le \ell <\Delta_i}
    \omega(w_\ell, \, w_{\ell-1}) \ge \ee^{- n^\varepsilon} \Big\}
 &=& [\Q \{ \omega(w_1, \, \varnothing) \ge \ee^{-n^\varepsilon} \}
    ]^{\Delta_i-1}
    \\
 &\ge&(1- 2\ee^{-n^\varepsilon})^{\Delta_i-1}
    \\
 &\ge& (1- 2\ee^{-n^\varepsilon})^n
    \ge
    1- \ee^{-n^{\varepsilon/2}}.
\end{eqnarray*}

\noindent We apply Lemma \ref{L:inter} to $A_i^{(n)} := \{ \omega(w_i, \, w_{i-1}) \ge \ee^{- n^\varepsilon}\}$ (with $\Delta_i-1$ and $(a_*-2\varepsilon)(n-m_i)^{1/3}$ playing the roles of $n$ and $r$, respectively; noting that $\Delta_i-1 \sim 3 (k-i)^2 n^\chi$ and $(n-m_i)^{1/3} \sim (k-i) n^{\chi/3}$, so the last condition in the lemma on $A_i^{(n)}$ is satisfied), to see that for all large $n$ and for all $1\le i\le k-1$,
$$
\Q\Big( F_i(n), \; \min_{1\le \ell <\Delta_i} \omega(w_\ell, \, w_{\ell-1}) \ge \ee^{- n^\varepsilon} \Big) \ge \exp\Big( - (1+\varepsilon ) {3\pi^2 \sigma^2 \over 8(a_* -2\varepsilon)^2} n^{\chi/3}\Big),
$$

\noindent which implies that
$$
\Q\Big( E_n^{(2)} \cap
    \bigcap_{i=1}^{k-1} E_{n,i}^{(3)}\Big)
\ge {1\over 2^k}\,  \exp\Big( - (k-1)(1+\varepsilon ) {3\pi^2 \sigma^2 \over 8(a_* -2\varepsilon)^2} n^{\chi/3}\Big).
$$

\noindent By definition, $k = \lfloor n^{(1-\chi)/3 }\rfloor$; hence
$$
\liminf_{n\to \infty} {1\over n^{1/3}} \log \Q\Big( E_n^{(2)} \cap
    \bigcap_{i=1}^{k-1} E_{n,i}^{(3)}\Big)
\ge
- (1+\varepsilon) {3\pi^2 \sigma^2 \over 8(a_* -2\varepsilon)^2}.
$$

\noindent This, together with (\ref{fineho3}), yields
$$
\liminf_{n \to\infty} \, {1\over n^{1/3}} \log \E(\varrho_n)
\ge
- \Big( {3\pi^2\sigma^2\over 8}\Big)^{1/3} .
$$

\noindent Since $\beta_n \ge \varrho_n$ (Lemma \ref{l:beta-rho}), we obtain (\ref{low-epn}), thus the lower bound in Theorem \ref{t:main}.

\appendix
\section{Appendix. Proofs of Lemmas \ref{markov-zj} and \ref{L:inter}}

\noindent {\it Proof of Lemma \ref{markov-zj}.} Although the lemma is deterministic, our proof is probabilistic. We note that the value of $b_1$ plays no role in the lemma.

Let $(\eta_i)_{i\ge0}$ be a Markov chain on $\{0, 1, \cdots, n\}$ with transition probabilities
$$
\P\Big( \eta_{i+1}= k \, \big| \, \eta_i=j\Big)=
\left\{
 \begin{array}{ll}
 {b_{j+1}\over 1+ b_{j+1}}  , & \hbox{if $k=j+1$,}
 \\
 \\
 {1\over 1+ b_{j+1}}, & \hbox{if $k=j-1$.}
 \end{array}
\right.
\qquad 0< j<n .
$$

\noindent [The transition probabilities from $j=0$ and $j=n$, having no importance, can be anything. For example, we can take $\P \{ \eta_{i+1}= 0 \, | \, \eta_i=0\} =1 = \P \{ \eta_{i+1}= n \, | \, \eta_i=n\}$.] Define $\tau_\eta(j):=\inf\{i\ge1: \eta_i=j\}$ (with $\inf \varnothing := \infty$ as usual). Let
$\P_j$ be the probability such that $\P_j\{\eta_0=j\} =1$, and let $\E_j$ be the expectation with respect to $\P_j$. We claim that ($\infty$ being not $<\infty$)
\begin{equation}
    \prod_{j=1}^{n-1} z_j
    \ge
    \prod_{i=1}^k \E_{m_i}
    \Big( {\bf 1}_{(\tau_\eta(m_{i-1}) < \tau_\eta(m_i))}
    (1+ r^{(i)})^{- \tau_\eta(m_{i-1})}\Big),
    \label{zjemi}
\end{equation}

\noindent and for any integers $0\le \ell<  m\le n$ and $r\ge 0$,
\begin{eqnarray}
    \E_m \Big( (1+r)^{-\tau_\eta(\ell)}
    {\bf 1}_{( \tau_\eta(\ell) < \tau_\eta(m))}\Big)
 &\ge&{1\over 2(m-\ell)}
    \exp\Big\{ - \max_{\ell<i\le m} (v(i) - v(m))
    \nonumber
    \\
 && \hskip-60pt
    - (v(m) - v(m+1))^+
    - r (m-\ell)^2 \ee^{\max_{\ell <i\le j\le m} [v(i)-v(j)]} \Big\}.
    \label{ablen}
\end{eqnarray}

\noindent Plainly Lemma \ref{markov-zj} will follow from (\ref{zjemi}) and (\ref{ablen}).

To prove (\ref{zjemi}), we consider a Markov chain $(\widetilde \eta_i)_{i\ge0}$ on $\{0, 1, \cdots, n\} \cup \{ \partial\}$, where $\partial$ is an absorbing point, such that
$$
\P\Big( \widetilde \eta_{i+1}= k \, \big| \, \widetilde\eta_i=j\Big)
=
\left\{
 \begin{array}{ll}
 b_{j+1}\, q_j , & \hbox{if $k=j+1$,}
 \\
 q_j, & \hbox{if $k=j-1$,}
 \\
 r_j q_j, & \hbox{if $k=\partial$,}
 \end{array}
\right. \qquad 0< j<n,
$$

\noindent with $q_j := {1\over b_{j+1} +1+r_j}$ for all $0<j < n$. [Again, the transition probabilities from $j=0$ and $j=n$ can be anything.] Let $\tau_{\widetilde\eta}(j):= \inf\{ i\ge 1: \widetilde\eta_i=j\}$. Then 
\begin{equation}
    z_j
    =
    \P_j \Big( \tau_{\widetilde\eta}(j-1)<
\tau_{\widetilde\eta}(n) \Big)  , 
    \qquad \forall 1\le j \le n-1.
    \label{A3}
\end{equation}

\noindent Let us check (\ref{A3}): $z_{n-1}= q_{n-1}$, and if $z_j= \P_j
(\tau_{\widetilde\eta}(j-1)< \tau_{\widetilde\eta}(n))$ for
$j\in [2, \, n-1]\cap \z$, then $\P_{j-1} ( \tau_{\widetilde\eta}(j-2)<
\tau_{\widetilde\eta}(n) ) = q_{j-1} + b_j q_{j-1} \P_j (
\tau_{\widetilde\eta}(j-2)< \tau_{\widetilde\eta}(n)) = q_{j-1}
+ b_j q_{j-1} z_j \P_{j-1} ( \tau_{\widetilde\eta}(j-2)<
\tau_{\widetilde\eta}(n) )$ by the Markov property. Hence $\P_{j-1} (
\tau_{\widetilde\eta}(j-2)< \tau_{\widetilde\eta}(n)) = {
q_{j-1} \over 1- b_j q_{j-1} z_j}$, which, by the definition of $q_{j-1}$, is ${1\over 1+ r_{j-1} +  b_j (1-z_j)
}$. Hence $\P_{j-1} ( \tau_{\widetilde\eta}(j-2)<
\tau_{\widetilde\eta}(n) )=z_{j-1}$. This establishes (\ref{A3}). As a consequence,
\begin{equation}
    \P_{n-1}\Big( \tau_{\widetilde\eta}(0)<
    \tau_{\widetilde\eta}(n)\Big)
    =
    \prod_{j=1}^{n-1} z_j.
    \label{pnt}
\end{equation}

We claim that for $0\le \ell < m <n$,
\begin{equation}
    \P_m \Big( \tau_{\widetilde\eta}(\ell)
     < \tau_{\widetilde\eta}(m)\Big)
    \ge
    \E_m \Big( {\bf 1}_{( \tau_{ \eta}(\ell) < \tau_{ \eta}(m))}
    (1+r)^{ -  \tau_{ \eta}(\ell) }\Big),
    \label{pkt}
\end{equation}

\noindent where $r:= \max_{\ell< j \le m} r_j$. Since $\P_{n-1} ( \tau_{\widetilde\eta}(0)< \tau_{\widetilde\eta}(n)) \ge \prod_{i=1}^k \P_{m_i} ( \tau_{\widetilde\eta}(m_{i-1})< \tau_{\widetilde\eta}(m_i ) )$, (\ref{zjemi}) will be a consequence of (\ref{pnt}) and (\ref{pkt}).

To prove (\ref{pkt}), let $\Omega_{\ell, m}$ be the set of all (finite)
paths of $\widetilde \eta$ starting from $m$ and hitting $\ell$ before
returning to $m$ (and without being absorbed by $\partial$). For any
$\gamma \in \Omega_{\ell, m}$, let $L_\gamma^\pm(j):= \sum_{i\ge
0} {\bf 1}_{(\gamma_i= j, \, \gamma_{i+1}= j\pm 1)}$ and $L_\gamma(j
) := L_\gamma^+(j) + L_\gamma^-(j)$. Then
\begin{eqnarray*}
    \P_m \Big( \tau_{\widetilde\eta}(\ell)
     < \tau_{\widetilde\eta}(m)\Big)
 &=& \sum_{\gamma \in \Omega_{\ell, m}}
    \prod_{\ell<j\le m} (b_{j+1} q_j)^{L_\gamma^+(j)}
    (q_j)^{L_\gamma^-(j)}
    \\
 &=& \sum_{\gamma \in \Omega_{\ell, m}}
    \prod_{\ell<j\le m}
    ({b_{j+1} \over 1+b_{j+1}})^{L_\gamma^+(j)}
    ({1\over 1+b_{j+1}})^{L_\gamma^-(j)}
    \Big( {1+b_{j+1} \over 1+b_{j+1} +r_j}\Big)^{L_\gamma(j)}
    \\
 &\ge& \sum_{\gamma \in \Omega_{\ell, m}}
    \prod_{\ell<j\le m}
    ({b_{j+1} \over 1+b_{j+1}})^{L_\gamma^+(j)}
    ({1\over 1+b_{j+1}})^{L_\gamma^-(j)}
    \Big( { 1  \over 1 + r }\Big)^{L_\gamma(j)}
    \\
 &=& \E_m \Big( {\bf 1}_{( \tau_{\eta}(\ell) < \tau_{ \eta}(m))}
    (1+r)^{ -  \tau_{ \eta}(\ell) }\Big),
\end{eqnarray*}

\noindent yielding (\ref{pkt}) and hence (\ref{zjemi}).

It remains to show (\ref{ablen}). We note that $\P_m \{ \tau_\eta(\ell) < \tau_\eta(m) \} = {1\over 1+ b_{m+1}} \P_{m-1} \{ \tau_\eta(\ell) < \tau_\eta( m) \}$. To compute $\P_{m-1} \{ \tau_\eta(\ell) < \tau_\eta( m) \}$ for the birth-and-death chain, we recall $v(j)= - \sum_{i=1}^j
\log b_i$ for $j\ge1$, and use the same argument as in the second identity of (\ref{zeitouni}) to see that $\P_{m-1} \{ \tau_\eta(\ell) < \tau_\eta( m) \} = {\ee^{v(m)}\over \sum_{i=\ell+1}^m \ee^{v(i)}}$. Therefore,
\begin{eqnarray}
    \P_m \Big(\tau_\eta(\ell) < \tau_\eta(m) \Big)
 &=&{1\over 1+ b_{m+1}} \,
    {\ee^{v(m)}\over \sum_{i=\ell+1}^m \ee^{v(i)}}
    \nonumber
    \\
 &\ge& {1\over 2(m-\ell)}
    \exp\Big\{ - \max_{\ell<i\le m} (v(i) - v(m)) - (v(m) - v(m+1))^+
    \Big\}.
    \label{m-l}
\end{eqnarray}

\noindent Under $\P_m$ and  conditionally on  $\{ \tau_\eta(\ell) <
\tau_\eta(m)\}$, $\tau_\eta(\ell)$ is stochastically smaller than the
hitting time of $\ell$ by a Markov chain  with the same probability
transition as $\eta$ but reflecting on $m$. The expectation of the
latter hitting time was estimated by Golosov (\cite{G84}, p.~498,
(A.1)). Hence
$$
\E_m \Big( \tau_\eta(\ell) \, \big| \, \tau_\eta(\ell) <
\tau_\eta(m) \Big) \le (m-\ell)^2 \exp\Big( \max_{\ell < i \le j \le  m} (v(i) - v(j))\Big),
$$

\noindent which, by means of the elementary inequality $(1+r)^{-u} \ge \ee^{-ru}$ for $u\ge 0$ and Jensen's inequality, implies that
$$
\E_m \Big( (1+r)^{-\tau_\eta(\ell)} \, \big| \, \tau_\eta(\ell) <
\tau_\eta(m) \Big)
\ge
\exp\Big( - r (m-\ell)^2 \ee^{\max_{\ell < i \le j \le  m} (v(i) - v(j))}\Big).
$$

\noindent  This together with (\ref{m-l}) implies (\ref{ablen}), completing the proof of Lemma \ref{markov-zj}.\hfill$\Box$

\bigskip

\noindent {\it Proof of Lemma \ref{L:inter}.} We start with the proof of the lower bound. Let $a\ge 2$ be an integer whose value will be chosen later on. Let $\eta < {1\over a}$. Let $K:= \lfloor {n \over a  r^2}\rfloor$, and $n_i= i a r^2$ for $0\le i\le K-1$ and $n_K:= n$. Write $S^{\#}_n:= \max_{1\le i \le n}( {\overline{S}_i - S_i})$. It is clear that
\begin{eqnarray*}
 && \P\Big( S^{\#}_n< r, \, \overline{S}_n= S_n, \,
    \bigcap_{i=1}^n A_i^{(n)}\Big)
    \\
 &\ge& \P\Big(\forall 1\le j \le K, \,
    \max_{n_{j-1} < i \le n_j} (\overline{S}_i - S_i) <  r, \,
    \overline{S}_{n_j} - S_{n_j} <  \delta r, \,
    \overline{S}_n= S_n,  \bigcap_{i=1}^n A_i^{(n)}\Big).
\end{eqnarray*}

For any $1\le j \le K$, conditionally on $\sigma\{S_i, A^{(n)}_i,
1\le i\le n_{j-1}\}$ and on $\{\overline{S}_{n_{j-1}} - S_{n_{j-1}} =
x\}$, the reflecting random walk  $( \overline{S}_{i+n_{j-1}}
-S_{i+n_{j-1}}, \, 0\le i \le n_j-n_{j-1})$ has the same law as
$(\max\{ x, \, \overline{S}_i\} - S_i, \, 0\le i \le n_j-n_{j-1})$. Accordingly,
$$
\P\Big( S^{\#}_n < r, \, \overline{S}_n= S_n, \, \bigcap_{i=1}^n A_i^{(n)}\Big)
\ge
q_{n,r}^{ K-1} b_{n,r},
$$

\noindent where
\begin{eqnarray*}
    q_{n,r}
 &:=& \P\Big( S^{\#}_{a r^2} < (1-\delta) r, \,
    \overline{S}_{a r^2} - S_{ar^2} < \delta r, \,
    \overline{S}_{a r^2} > \delta r , \,
    \bigcap_{i=1}^{a r^2} A_i^{(n)} \Big) ,
    \\
    b_{n,r}
 &:=& \P\Big(S^{\#}_{n_K - n_{K-1}} < (1-\delta) r,\,
    \overline{S}_{n_K - n_{K-1}} = S_{n_K - n_{K-1}} > \delta r ,\,
    \bigcap_{i=1}^{n_K- n_{K-1}} A_i^{(n)}\Big).
\end{eqnarray*}

\noindent We observe that
\begin{eqnarray*}
    q_{n, r}
 &\ge& \P\Big( S^{\#}_{a r^2} < (1- \delta) r, \,
    \overline{S}_{a r^2} - S_{ar^2} < \delta r, \,
    \overline{S}_{a r^2} > \delta r  \Big)
    +
    \P\Big( \bigcap_{i=1}^{a r^2} A_i^{(n)} \Big)- 1
    \\
 &\ge& \P\Big( S^{\#}_{a r^2} < (1-\delta) r, \,
    \overline{S}_{a r^2} - S_{a r^2} < \delta r, \,
    \overline{S}_{a r^2} > \delta r  \Big) - {\eta \over r}.
\end{eqnarray*}

\noindent On the other hand, since the three events $\{S^{\#}_{a r^2} <
(1-\delta) r\}$, $\{\overline{S}_{a r^2} - S_{ar^2} < \delta r\}$ and $\{ \overline{S}_{a r^2} > \delta r\}$ are non-decreasing with respect to each $S_i- S_{i-1}$ (for $1\le i \le n$), it follows from the FKG inequality that
$$
q_{n, r}
\ge
\P\Big( S^{\#}_{a r^2} < (1-\delta) r \Big) \, \P\Big( \overline{S}_{a r^2} - S_{a r^2} < \delta r \Big) \, \P\Big( \overline{S}_{a r^2} > \delta r\Big) - {\eta \over r}.
$$

\noindent By Donsker's invariance principle,
\begin{eqnarray*}
 && \lim_{r\to\infty}
    \P\Big( S^{\#}_{a r^2} < (1-\delta) r \Big) \,
    \P\Big( \overline{S}_{a r^2} - S_{a r^2} < \delta r \Big) \,
    \P\Big( \overline{S}_{a r^2} > \delta r\Big)
    \\
 &=& \P\Big( \sup_{t\in [0, \, 1]} (\overline{W}(t) - W(t) )
    < {1-\delta \over \sigma \sqrt{a}} \Big) \,
    \P\Big(\overline{W}(1) -W(1) < {\delta\over \sigma\sqrt{a}} \Big)\,
    \P\Big( \overline{W}(1) > {\delta\over \sigma \sqrt{a}} \Big),
\end{eqnarray*}

\noindent  where $W$ is a standard Brownian motion, with $\overline{W}(t) := \sup_{s\in [0, \, t]} W(s)$ for any $t$. Recall that $\sup_{t\in [0, \, 1]} (\overline{W}(t) - W(t) )$ has the same distribution as $\sup_{t\in [0, \, 1]} |W(t)|$; so (see Formula (5.9), page 342 of Feller~\cite{feller}) $\P \{ \sup_{t\in [0, \, 1]} (\overline{W}(t) - W(t) ) \le x \}= \exp\{ - (1+o(1)) {\pi^2 \over 8 x^2}\}$, for $x\to0$. Consequently, for all sufficiently large $a$, say $a \ge a_0$, we have
$$
\P\Big( \sup_{t\in [0, \, 1]} (\overline{W}(t) - W(t) )  < {1-\delta \over \sigma \sqrt{a}} \Big)
\ge
\ee^{- (1+3\delta) {\pi^2 \sigma^2 a \over 8}}.
$$

\noindent Since both $\overline{W}(1) - W(1)$ and $\overline{W}(1)$
are distributed as the absolute value of a standard Gaussian random
variable, we can even increase the value of $a_0$ (if necessary) such that for all $a\ge a_0$ (thus ${\delta \over \sigma\sqrt{a}}$ is very small), $\P\{ \overline{W}(1) > {\delta \over \sigma \sqrt{a}} \} \ge {1\over 2}$ and $\P\{\overline{W}(1) - W(1) < {\delta \over \sigma\sqrt{a}} \} \ge {\delta\over 2 \sigma \sqrt {a}}$.

Now we fix an arbitrary integer $a\ge a_0$. For all large $r$ (say $r\ge r_0$) such that $r \le \eta \sqrt{n}$, we have
$$
q_{n, r} \ge  \ee^{- (1+4\delta) {\pi^2 \sigma^2 a \over 8}}.
$$

The probability $b_{n,r}$ can be estimated in a similar way: from the assumptions on $(A_i^{(n)})$ and the FKG inequality, we deduce that  $$
b_{n , r}
\ge
\P (S^{\#}_{n_K - n_{K-1}} < (1-
\delta) r ) \, \P( \overline{S}_{n_K - n_{K-1}} = S_{n_K - n_{K-1}} ) \, \P ( \overline{S}_{n_K - n_{K-1}} > \delta
r   )  - {\eta \over r} .
$$

\noindent Observe that $ar^2 \le n_K - n_{K-1} \le 2 a r^2$. Therefore, by Donsker's invariance principle, for all $r \ge r_0$ (with an increased value of $r_0$ if necessary), $\P (S^{\#}_{n_K - n_{K-1}} < (1- \delta) r ) \, \P( \overline{S}_{n_K - n_{K-1}} > \delta r) \ge c(a, \delta)$ for some constant $c(a, \delta)>0$, whereas $\P( \overline{S}_{n_K - n_{K-1}} = S_{n_K - n_{K-1}} )= \P( S_1\ge 0, S_2 \ge 0, ..., S_{n_K - n_{K-1} -1} \ge 0)  \ge {c' \over r\sqrt{a}}$ for some constant $c'>0$. Taking $\eta := \min\{ {c(a, \delta) c' \over 2 \sqrt{a}}, \, {1\over a }\}$, we get $ b_{n , r} \ge {\eta \over r}$. Consequently,
$$
\P \Big(S^{\#}_n < r, \, \overline{S}_n= S_n, \, \bigcap_{i=1}^n
A_i^{(n)} \Big)
\ge
q_{n,r}^{ K-1} b_{n,r}
\ge
{\eta \over r } \ee^{- (1+4\delta) { \pi^2 \sigma^2 \over 8} {n \over r^2} } ,
$$

\noindent proving the lower bound in the lemma.

The upper bound is easier: with the same notation and the choice of
$a$, we have
\begin{eqnarray*}
    \P\Big(S^{\#}_n < r, \, \overline{S}_n= S_n, \,
    \bigcap_{i=1}^n A_i^{(n)}\Big)
 &\le& \P\Big( \forall 1\le j \le K-1 , \,
    \max_{n_{j-1} < i \le n_j} (\overline{S}_i - S_i) <  r \Big)
    \\
 &\le& \Big[ \P (S^{\#}_{a r^2} < r)\Big]^{K-1} .
\end{eqnarray*}

\noindent For all large $r$, $\P(S^{\#}_{a r^2} < r) \le \ee^{- (1-\delta) {\pi^2 \sigma^2 a \over 8}}$; hence $\P (S^{\#}_n < r, \, \overline{S}_n= S_n, \, \cap_{i=1}^n A_i^{(n)} ) \le \ee^{- (1-\delta) {\pi^2 \sigma^2  \over 8 } { {n \over r^2 } -1}}$, yielding the upper bound (by eventually modifying the choice of $\eta$ in terms of $a$). \hfill$\Box$

\bigskip
\bigskip

\noindent {\Large\bf Acknowledgement}

\bigskip

\noindent We have been kindly informed by Ofer Zeitouni that Theorem \ref{t:max-ray} in the case of regular trees was independently proved by Fang and Zeitouni \cite{FZ09}. We are also grateful to an anonymous referee for her/his very careful reading of the original manuscript.

\bigskip
\bigskip


{\footnotesize

\baselineskip=12pt

\noindent
\begin{tabular}{lll}
& \hskip20pt Gabriel Faraud
    & \hskip60pt Yueyun Hu \\
& \hskip20pt D\'epartement de Math\'ematiques
    & \hskip60pt D\'epartement de Math\'ematiques \\
& \hskip20pt Universit\'e Paris XIII
    & \hskip60pt Universit\'e Paris XIII \\
& \hskip20pt 99 avenue J-B Cl\'ement
    & \hskip60pt 99 avenue J-B Cl\'ement \\
& \hskip20pt F-93430 Villetaneuse
    & \hskip60pt F-93430 Villetaneuse \\
& \hskip20pt France
    & \hskip60pt France \\
& \hskip20pt {\tt faraud@math.univ-paris13.fr}
    & \hskip60pt
    {\tt yueyun@math.univ-paris13.fr}
\end{tabular}

\bigskip
\bigskip

\hskip20pt Zhan Shi\par \hskip20pt Laboratoire de Probabilit\'es UMR
7599\par \hskip20pt Universit\'e Paris VI\par \hskip20pt 4 place
Jussieu\par \hskip20pt F-75252 Paris Cedex 05\par \hskip20pt
France\par \hskip20pt {\tt zhan.shi@upmc.fr}

}

\end{document}